\documentclass{amsart}
\usepackage[T1]{fontenc}
\usepackage{graphicx}
\usepackage{amsmath}
\usepackage{amssymb}
\usepackage{mathrsfs}
\usepackage{amsthm}
\usepackage{thmtools}
\usepackage{hyperref}
\usepackage{cleveref}
\usepackage{geometry}
\usepackage{tikz}
\usepackage{tikz-cd}
\usepackage{todonotes}
\usepackage{cite}
\usepackage{booktabs}
\usepackage{colortbl}
\ifcsname theorem\endcsname{}\else\declaretheorem[parent=section]{theorem}\fi
\ifcsname corollary\endcsname{}\else\declaretheorem[sibling=theorem]{corollary}\fi
\ifcsname lemma\endcsname{}\else\declaretheorem[sibling=theorem]{lemma}\fi
\ifcsname proposition\endcsname{}\else\declaretheorem[sibling=theorem]{proposition}\fi
\ifcsname conjecture\endcsname{}\else\fi
\ifcsname problem\endcsname{}\else\fi
\ifcsname question\endcsname{}\else\fi
\ifcsname definition\endcsname{}\else\fi
\ifcsname exercise\endcsname{}\else\fi
\ifcsname example\endcsname{}\else\declaretheorem[sibling=theorem, style=definition]{example}\fi
\ifcsname remark\endcsname{}\else\declaretheorem[sibling=theorem, style=remark]{remark}\fi

\providecommand {\Z}{{\bf Z}}
\providecommand {\Q}{{\bf Q}}
\providecommand {\R}{{\bf R}}
\providecommand {\C}{{\bf C}}

\renewcommand {\P}{{\bf P}}

\providecommand {\A}{{\bf A}}

\providecommand{\GL}{\operatorname{GL}}
\providecommand{\PGL}{\operatorname{PGL}}
\providecommand{\Gm}{{\bf G}_m}

\providecommand{\spec}{\operatorname{Spec}}

\providecommand{\Hom}{\operatorname{Hom}}

\providecommand{\codim}{\operatorname{codim}}
\providecommand{\Sym}{\operatorname{Sym}}

\renewcommand{\k}{\mathbf k}
\DeclareMathOperator{\G}{\mathbf G}
\DeclareMathOperator{\Orb}{Orb}
\DeclareMathOperator{\sP}{\mathscr P}
\DeclareMathOperator{\Bl}{Bl}
\renewcommand{\O}{\mathcal O}

\date{\today}
\begin{document}

\title{Equivariant classes of orbits in \(\GL(2)\)-representations}
\author{Anand Deopurkar}
\address{Mathematical Sciences Institute, Australian National University, Canberra, Australia}

\begin{abstract}
  We compute equivariant fundamental classes of orbits in \(\GL(2)\)-representations.
  As applications, we find degrees of the orbit closures corresponding to elliptic fibrations and self-maps of the projective line.
\end{abstract}

\maketitle

\section{Introduction}
If we fix a hypersurface in projective space, how complicated is the set of all hypersurfaces obtained from the fixed one by changes of coordinates?
Similarly, if we fix a self-map of the projective space, how complicated is the set of all self-maps obtained from the fixed one by changes of coordinates?
These questions, and many others, generalise as follows.
Given a representation \(W\) of an algebraic group \(G\), how complicated is the \(G\)-orbit of a fixed \(w \in W\)?
One measure of complexity is the degree of the orbit closure in \(\P W\).
A more refined measure is the \(G\)-equivariant fundamental class.
Our main theorem (\Cref{thm:main}) completely describes the equivariant fundamental classes (and hence degrees) of orbits in representations of \(G = \GL(2)\).
The case where \(W\) is irreducible was already known.
The new contribution is treating reducible \(W\); this presents new challenges, but also has new applications.
The case of \(G = \GL(1)\), or more generally, any torus, is straightforward.
We treat it in \Cref{sec:torus}.

The question of finding equivariant classes of orbit closures has been well studied, especially in cases where the orbits have a geometric interpretation.
For the \(\GL(2)\) representation \(\Sym^n \C^2\), where the orbits represent divisors of degree \(n\) on \(\P^1\) modulo changes of coordinates, the degree of the orbit closure was computed by Enriques--Fano \cite{enr.fan:98} for the generic case and Aluffi--Faber in general \cite{alu.fab:93*1}.
The equivariant class was computed by Lee--Patel--Tseng \cite[Appendix~B]{lee.pat.tse:23}.
For the \(\GL(3)\) representation \(\Sym^n \C^3\), where the orbits represent plane curves of degree \(n\), the degree of the orbit closure was computed by Aluffi--Faber \cite{alu.fab:00,alu.fab:93}.
For the \(\GL(4)\) representation \(\Sym^3 \C^4\), where the orbits represent cubic surfaces, the equivariant class of a generic orbit closure was computed in \cite{deo.pat.tse:21}.
Local analogues of equivariant orbit classes are Thom polynomials, which have been studied by Buch, Fehér, Rimányi, and Weber among others \cite{feh.rim.web:20,feh.rim:12,rim:01}.
In all these cases, the equivariant class yields a counting formula---the equivariant orbit/Thom class of \(w\) gives the number of times \(w\) appears, up to isomorphism, in a given family.
We expect the equivariant class to reflect the geometry of \(w\).
This is indeed the case for divisors on \(\P^1\), where the class depends on the multiplicities in the divisor, and for curves in \(\P^2\), where the class depends on the singularities and flexes of the curve.

As a direct application of the main theorem, we compute the degrees of orbit closures in two (reducible) representations of geometric significance.
The first is the \(\GL(2)\)-representation \(\Sym^{4n}(\C^2) \oplus \Sym^{6n}(\C^2)\), where the orbits represent isomorphism classes of elliptic fibrations over \(\P^1\).
In this case, the degree depends on the Kodaira types of the singular fibers.
The second is the \(\GL(2)\)-representation \(\Hom(\C^2, \Sym^n\C^2)\), where most orbits represent isomorphism classes of self-maps \(\P^1 \to \P^1\) of degree \(n\).
(The space of degree \(n\) maps \(\P^1 \to \P^1\) and its quotient by changes of coordinates are important objects of study in complex dynamics, where they are usually denoted by \(\operatorname{Rat}_n\) and \(\mathcal M_n\) \cite{ree:92,seg:79,mil:93,sil:98}).
In this case, the degree is (surprisingly) independent of the orbit.

We first give these two applications in \Cref{sec:introellipticfibrations} and \Cref{sec:introrationalmaps}, respectively, before stating the main theorem in \Cref{sec:intromain}.

\subsection{Elliptic fibrations}\label{sec:introellipticfibrations}
Fix a positive integer \(n\), and let \(W = \Sym^{4n}(\C^2) \oplus \Sym^{6n}(\C^2)\).
A non-zero \((A, B) \in W\) determines an elliptic fibration over \(\P^{1}\) defined by the Weierstrass equation
\[ y^2 = x^3 + Ax + B.\]
The \(\GL(2)\)-orbits in \(W\) thus represent Weierstrass elliptic fibrations over \(\P^1\), up to isomorphism.

Let \(h\) be the class of the Weil divisor \(\O(1)\) on the weighted projective space \(\P W = (W - 0) \big / \Gm\), where \(\Gm\) acts by weight \(2\) on \(\Sym^{4n}(\C^2)\) and by weight \(3\) on \(\Sym^{6n}(\C^2)\).
Given \((A, B) \in W\) and \(u \in \P^{1}\), let \(\operatorname{ord}(A)_u\) and \(\operatorname{ord}(B)_u\) be the orders of vanishing of \(A\) and \(B\) at \(u\).
Set
\[c(u) = \min\left(\frac{1}{2}\operatorname{ord}(A)_u, \frac{1}{3}\operatorname{ord}(B)_u\right).\]

\begin{theorem}\label{thm:ellipticfibrations}
  Fix a non-zero \(w = (A,B) \in  W = \Sym^{4n}(\C^2) \oplus \Sym^{6n}(\C^2)\), and let \(\pi \colon E \to \P^1\) be the Weierstrass fibration defined by \(w\).
  Let \(D \subset \P^1\) be a finite set such that \(\pi\) is smooth on \(\P^{1} - D\).
  Let \(\Gamma \subset \GL(2)\) be the stabiliser of \(w \in W\) and let \(\Orb([w])\) be the closure of the \(\PGL(2)\)-orbit of \([w] \in \P W\).
  Then
  \[ |\Gamma| [\Orb([w])] = 2^{4n+3}3^{6n+1}\cdot n \cdot \left(4n^3-\sum_{u \in D}c(u)^2(3n-c(u)) \right) h^{10n-2}.\]
\end{theorem}
If \(\pi\) is a minimal Weierstrass fibration as in \cite[III.3]{mir:89}, then \(c(u) < 2\) and \(c(u)\) determines the Kodaira fiber type over \(u\) (see \cite[IV.3.1]{mir:89}).
See \Cref{tab:ellipticfibers} for the Kodaira types and their contribution to the formula above.
The main theorem in fact gives the equivariant class, of which the degree is a particular specialisation.
See \Cref{sec:ellipticfibrations} for the proof.
\begin{table}[h]
  \centering
  \rowcolors{2}{white}{gray!15}
  \begin{tabular}[h3]{llp{.42\textwidth}p{.3\textwidth}}
    \toprule
\(c(u)\) &    Type & Description  & Contribution to the degree \(c(u)^2(3n-c(u))\)\\
    \midrule
\(0\)&    \(I_N\) &Smooth elliptic curve, nodal rational curve, or cycle of smooth rational curves& \(0\) \\
\(1\) &    \(I_N^{*}\)& \(\widetilde D_{4+N}\)-configuration of rational curves &\(3n-1\) \\
\(1/3\) &    \(II\) &Cuspidal rational curve & \({1}/{27} \cdot (9n-1)\)\\
\(1/2\) &    \(III\)&Two tangent rational curves&\({1}/{8} \cdot (6n-1)\) \\
\(2/3\) &    \(IV\)&Three concurrent rational curves &\({4}/{27} \cdot (9n-2)\) \\
\(4/3\)&    \(IV^{*}\)&\(\widetilde E_{6}\)-configuration of rational curves & \({16}/{27} \cdot (9n-4)\)\\
\(3/2\)&    \(III^{*}\)&\(\widetilde E_{7}\)-configuration of rational curves &\({27}/{8} \cdot (2n-1)\) \\
\(5/3\)&    \(II^{*}\)&\(\widetilde E_{8}\)-configuration of rational curves &\({25}/{27} \cdot (9n-5)\)\\
    \bottomrule
  \end{tabular}
    \caption[Contributions of singular fibers]{Contributions from the singular fibers in a minimal Weierstrass fibration \(y^2 = x^3 + Ax +B\) towards the degree of the orbit closure of \((A,B) \in \P\left(\Sym^{4n}(\C^{2}) \oplus \Sym^{6n}(\C^{2})\right)\). }
  \label{tab:ellipticfibers}
\end{table}

\subsection{Rational self maps}\label{sec:introrationalmaps}
Fix a positive integer \(n\) and set \(W = \Hom(\C^2, \Sym^n\C^2)\).
An element \(f \in W\) is equivalent to a map
\begin{equation}\label{eqn:endrat}
  \C^{2} \otimes \O_{\P^1} \to \O_{\P^{1}}(n).
\end{equation}
For \(f\) in a Zariski open subset, the map \eqref{eqn:endrat} is surjective, and hence defines a map \(\P^{1}  \to \P^{1}\) of degree \(n\).
Conversely, every map \(\P^{1} \to \P^{1}\) of degree \(n\) arises from an \(f \in W\), which is unique up to a scalar.
Thus, most \(\GL(2)\)-orbits in \(W\) represent maps \(\P^1 \to \P^1\) of degree \(n\) modulo changes of coordinates.

\begin{theorem}\label{thm:ratmaps}
  Suppose \(f \in \Hom(\C^2, \Sym^n\C^2)\) defines a map \(\P^{1}\to \P^{1}\) of degree \(n\).
  Let \(\overline \Gamma \subset \PGL(2)\) be the stabiliser of \([f] \in \P \Hom(\C^2, \Sym^n\C^2)\), and \(\Orb([f])\) the closure of the \(\PGL(2)\)-orbit of \([f]\).
  Then
  \[ |\overline \Gamma| \cdot \deg(\Orb([f])) = n(n+1)(n-1).\]
\end{theorem}
Again, the main theorem gives the equivariant class, of which the degree is a particular specialisation.
See \Cref{sec:ratmaps} for the proof.

We highlight that \Cref{thm:ratmaps} does not hold for all \(f \in \Hom(\C^2, \Sym^n\C^2)\).
  For the \(f\) whose associated map \eqref{eqn:endrat} is not surjective---\(f\) with base-points---the stabiliser-weighted degree can be different.
  It is remarkable that for the \(f\) without base-points, it is constant.
  This is in contrast to the case of divisors on \(\P^1\), where the multiplicities in the divisor matter.

\subsection{Main theorem}\label{sec:intromain}
Fix a 2-dimensional vector space \(V\) over an algebraically closed  field \(\k\) of characteristic 0.
Fix a finite dimensional \(\GL V\) representation
\[ W = W_{1} \oplus \cdots \oplus W_{n}, \text{ where } W_i = \Sym^{a_i-b_i}V \otimes \det{V}^{b_i}.\]
Set \(d_i = a_i+b_i\), and assume that \(d_i > 0\) for all \(i\).
Fix a maximal torus \(T \subset \GL V\).
We then have an isomorphism between the equivariant (rational) Chow ring \(A_{\GL V}\) and the symmetric polynomials in \(A_T = \Q[v_1,v_2]\).

We must now introduce some notation.
Fix a non-zero \(w = (w_1, \dots, w_n) \in W\), and write \(w_i = f_i \otimes \delta^{b_i}\) for some \(f_i \in \Sym^{a_i-b_i}V\) and \(\delta \in \det V\).
Given \(u \in \P^1\), let \(r_{i}^u\) be the order of vanishing of \(f_i\) at \(u\).
Let \(\Lambda^u \subset \R^2\) be the convex hull of the union of the shifted quadrants
\[ \frac{1}{d_i} (r_i^u+b_i,b_i) + \R_{\geq 0}^2.\]
Let \(\lambda^u(0), \dots, \lambda^u(k^u)\) be the vertices of \(\Lambda^u\) arranged from the bottom right to the top left.
For a \(p \in \R^2\), use \(p_{1}\) and \(p_2\) to denote the first and the second coordinates.
Set
\begin{align*}
  b = \min\left(b_i/d_i \mid w_i \neq 0\right), \qquad r_{\rm gen}^u &= \lambda^u(0)_1 - b, \qquad r^u = \min\left((r_i^u+b_i)/d_i\right), \text{ and } \\
  s^u &= \begin{cases}
    1-\frac{\lambda^u(0)_1 - \lambda^u(1)_1}{\lambda^u(0)_2 - \lambda^u(1)_2} \text{, if \(k^u \geq 1\),} \\
    1 \text{, otherwise}.
    \end{cases}
\end{align*}
For \(j = 1, \dots, k^u\), let \(\eta^u(j)\) and \(\zeta^u(j)\) be the smallest integral normal vectors to the rays of \(\Lambda^u\) at the vertex \(\lambda^u(j)\).
Set \(N^u(j) = \det(\eta^u(j), \zeta^u(j))\).
Let \(A \subset \P^1\) be a finite set that includes the common zero locus of \(\{f_i \mid b_i/d_i = b\}\).

Given \(F \in \Q(v_1,v_2)\), denote by \(F_{\rm sym}\) its symmetrisation
\[ F_{\rm sym} = F(v_1,v_2) + F(v_2,v_1).\]
Let \(N = \dim W\) and observe that in \(A_T = \Q[v_1,v_2]\), we have the top Chern class
\[ c_N(W) = \prod_{i=1}^n\prod_{j = 0}^{a_i-b_i}\left((b_{i}+j)v_1 + (a_i-j)v_2\right).\]
Let \(\Gamma \subset \GL(V)\) be the stabiliser of \(w \in W\); assume that it is finite.
\begin{theorem}\label{thm:main}
  In the notation above, the \(\GL V\)-equivariant class of the orbit closure of \(w \in W\) in \(A_{\GL V}(W) \subset \Q[v_1,v_2]\) is given by
  \begin{equation}\label{eqn:main}
    |\Gamma|[\Orb(w)] =    c_N(W) \cdot\left( F_{\rm sym} + \sum_{u \in A} G^u_{\rm sym} + \sum_{u \in A} \sum_{j = 1}^{k^u} H^u(j)_{\rm sym}\right),
  \end{equation}
  where
  \begin{align*}
    F &= 2((1-b)v_1+bv_2)^{-1}(v_1-v_2)^{-3} \\
    & \qquad - (2b-1)((1-b)v_1+v_2)^{-2}(v_1-v_2)^{-2} \\
    G^u &= ((1-r^u)v_1 + r^uv_2)^{-1}(v_1-v_2)^{-3} \\
      & \qquad -s^u((1-b)v_1+bv_2)^{-1}(v_1-v_2)^{-3} \\
    &\qquad -r_{\rm gen}^u((1-b)v_1+v_2)^{-2}(v_1-v_2)^{-2}, \text{ and }\\
    H^u(j) &= |N^u(j)|\eta^u(j)_1^{-1}\zeta^u(j)_{1}^{-1}((1-\lambda^u(j)_2)v_1 + \lambda^u(j)_2v_2)^{-1}(v_1-v_2)^{-3}.
  \end{align*}    
\end{theorem}
Note that in the sum of \(H^u(j)_{\rm sym}\), the bottom right vertex (\(j = 0\)) is omitted.
\begin{remark}
  It is not obvious that the expression in \Cref{thm:main} is a polynomial.
  But it must be, as a consequence of the theorem.
\end{remark}

\subsection{Negative or mixed weights}\label{sec:mixed}
Our main theorem applies to representations \(W\) whose direct summands have positive weights \(d_i\).
The theorem can also be used for \(W\) whose direct summands have negative weights by dualising or by twisting by a large negative \(n\) as described in \Cref{sec:twist}.

The cases where \(W\) has summands of weight 0 or some summands of positive weights and some of negative weights are a bit strange.
In these cases, a generic \(w \in W\) does not contain the origin in its orbit closure.
Therefore, its equivariant class of a generic orbit closure is \(0\), as can be seen by pulling back to the equivariant Chow ring of the origin.

\subsection{Ideas in the proof}
Let \(\P W\) be the weighted projective space \((W - 0)\big /\Gm\) for the central \(\Gm \subset \GL(2)\).
Given a \(w \in W\), the key idea is to find a complete orbit parametrisation for \(\Orb([w])\), namely a proper \(\PGL(2)\)-variety \(X\) and an equivariant finite map \(X \to \P W\) whose image is \(\Orb([w])\).
Then the class of \(\Orb([w])\) is the push-forward of \([X]\), up to a constant factor.
The push-forward also gives \(\GL(2)\)-equivariant class of \(\Orb(w)\) (see \Cref{prop:weightedformula}).

To find \(X\), we start with \(M = \P \Hom(\k^2, \k^2)\), and the rational map \(M \dashrightarrow \P W\) given by \(m \mapsto m w\).
We find an explicit resolution \(\widetilde M \to \P W\), which serves as our complete orbit parametrisation.
We then compute the push-forward as an integral on \(\widetilde M\) using Atiyah--Bott localisation.

The resolution \(\widetilde M \to M\) is a weighted blow-up.
It is much more convenient to take the weighted blow-up in a stacky sense.
The stacky blow-up is smooth and maps to the weighted projective stack \(\sP W\).
We can then write the push-forward as an integral and evaluate it using localisation.
The stacky blow-up is toroidal, and is completely described by the combinatorial data of the Newton polygons \(\Lambda^u\).

\subsection{Conventions and organisation}
We work over an algebraically closed field \(\k\) of characteristic 0.
A stack means an algebraic stack over \(\k\).
All schemes and stacks are of finite type over \(\k\).
Given a vector space/bundle \(V\), the projectivisation \(\P V\) refers to the space of one-dimensional sub-spaces/bundles of \(V\), consistent with the convention in \cite{ful:98} and \(\O_{\P V}(-1)\) denotes the universal sub-bundle.
All Chow groups are with rational coefficients.

In \Cref{sec:stackyblowup}, we recall stacky weighted blow-ups in preparation for our main construction.
In \Cref{sec:param}, we describe how to find the equivariant class of an orbit using a complete parametrisation.
Both of these sections are general (not specific to \(\GL(2)\)).
In \Cref{sec:twist}, we observe that the main theorem is invariant under a twist operation, which allows some simplification.
In \Cref{sec:completeparam}, we construct a complete parametrisation of a \(\GL(2)\)-orbit using a stacky blow-up.
In \Cref{sec:localisation}, we evaluate the equivariant orbit class using localisation.
In \Cref{sec:applications}, we deduce the applications to elliptic fibrations and rational self maps.
In \Cref{sec:torus}, we explain the case of \(G\) a torus.

\subsection*{Acknowledgements}
The project arose from conversations with Anand Patel, to whom I am deeply grateful.
I thank Ming Hao Quek for sharing his expertise on stacky blow-ups.
I was supported by the grant \texttt{DE180101360} from the Australian Research Council.

\section{Rational Newton polyhedra and weighted blow-ups}\label{sec:stackyblowup}
The material in this section should be well-known to experts (see, for example, \cite[\S~2]{que:24}).

Set \(M = \Z^n\) and \(N = \Hom(M,\Z).\)
Let \(M_{\geq 0}\) be the set of vectors with non-negative coordinates in \(M \otimes \R = \R^n\), and similarly for \(N_{\geq 0}\).
A \emph{rational Newton polyhedron} is a closed convex polyhedron \(\Lambda \subset M\) whose recession cone is \(M_{\geq 0}\) and whose vertices have rational coordinates.
Such a \(\Lambda\) gives a fan \(\Lambda^\perp\) in \(N \otimes \R\) supported on \(N_{\geq 0}\), called the \emph{normal fan} of \(\Gamma\).
There is an inclusion reversing bijection between the faces of \(\Lambda\) and the cones of \(\Lambda^\perp\).
To a face \(F\) of \(\Lambda\), we associate the cone \(F^{\perp}\) of \(\Lambda^{\perp}\) defined by
\[ F^{\perp}  = \{f \in N \mid f \text{ is constant on } F \text{ and this constant is the minimum of }f\text{ on } \Lambda\}.\]
Since the recession cone of \(\Lambda\) is \(M_{\geq 0}\), and \(f\) achieves a minimum on \(\Lambda\), it must lie in \(N_{\geq 0}\).

Let \(F\) be a maximal face of \(\Lambda\), that is, of dimension \((n-1)\).
Then \(F^{\perp}\) is a ray.
For every \(F\), choose a non-zero vector \(\beta_{F} \in F^{\perp}\) with integer coordinates.
Let \(r\) be the number of maximal faces of \(\Lambda\).
Then the collection \(\{\beta_F\}\) gives a homomorphism \(\beta \colon \Z^r \to N\) with finite cokernel.
Let \(\mathscr X_{\Lambda, \beta}\) be the toric stack defined by the data \((N, \Lambda^\perp, \beta)\) in the sense of \cite{bor.che.smi:05}.
It comes with a canonical map \(\mathscr X_{\Lambda, \beta} \to \A^n\), which we call the stacky blow-up of \(\A^n\) defined by \((\Lambda,\beta)\).

Let us describe \(\mathscr X_{\Lambda, \beta} \to \A^n\) in charts, following \cite[Proposition~4.3]{bor.che.smi:05}.
Assume that \(\Lambda\) is simplicial, that is, every vertex of \(\Lambda\) has exactly \(n\) incident rays.
Let \(v\) be a vertex of \(\Lambda\).
Denote the rays incident to \(v\) by \(R_1, \dots, R_n\) and the maximal faces incident to \(v\) by \(F_1, \dots, F_n\) such that \(R_i\) is the only ray not contained in \(F_i\).
Set \(\beta_i = \beta_{F_i}\) and let \(r_i \in R_i\) be the unique vector such that \(\langle \beta_i, r_i \rangle = 1\).
Then \(r_1, \dots, r_n\) is a basis of \(M \otimes \Q\) dual to the basis \(\beta_1, \dots, \beta_n\) of \(N \otimes \Q\).
Let \(M_v \supset M\) be the dual lattice of the sub-lattice of \(N\) spanned by \(\beta_1, \dots, \beta_n\).
Then \(M_v/M\) is a finite abelian group.
Set \(\mu_v = \Hom(M_v/M, \G_m)\).
The chart of \(\mathscr X_{\Lambda, \beta}\) defined by \(v\) is
\begin{equation}\label{eqn:chart}
  [\spec \k[u_1, \dots, u_n] / \mu_v],
\end{equation}
with the action given as follows.
A \(\zeta \in \mu_v\) acts by
\[\zeta \colon u_i \mapsto \zeta(r_i) u_i.\]
In particular, note that \(\mathscr X_{\Lambda, \beta}\) is a smooth Deligne--Mumford stack.

Let \(e_1, \dots, e_n\) be the standard basis vectors in \(M\).
In the chart above, the map to \(\A^n = \spec \k[x_1,\dots,x_n]\) is defined by
\begin{equation}\label{eqn:blowupmap}
  x_i \mapsto u_1^{\langle \beta_1, e_i \rangle} \cdots u_n^{\langle \beta_n, e_i \rangle}.
\end{equation}
Note that \(\zeta \in \mu_v\) multiplies the image of \(x_i\) by \(\zeta(e)\) where
\[ e = r_1 \langle \beta_1, e_i  \rangle + \dots + r_n\langle  \beta_n,e_i\rangle.\]
Since \(r_1, \dots, r_n\) and \(\beta_1,\dots,\beta_n\) are dual bases, we see that \(e = e_i \in M\) and hence \(\zeta(e) = 1\).
So the map \eqref{eqn:blowupmap} is indeed \(\mu_v\)-invariant.
Write \(r_i = (a_1,\dots,a_n)\) in standard coordinates with \(a_i \in \Q\).
Informally, it is helpful to think of 
\( u_i\) as \(x_1^{a_1} \cdots x_n^{a_n}\).

Let \(X_{\Lambda}\) be the toric variety associated to \((N, \Lambda^\perp)\).
Then we have a map
\(\mathscr X_{\Lambda, \beta} \to X_{\Lambda}\)
which is the coarse space map \cite[Proposition~3.7]{bor.che.smi:05}.

\begin{remark}\label{rem:eqvblowup}
  Let \((N, \Lambda^{\perp}, \beta)\) be the stacky fan given by a rational Newton polyhedron as above.
  Let \(r\) be the number of rays of \(\Lambda^{\perp}\).
  In \cite{bor.che.smi:05}, the stack associated to \((N, \Lambda^{\perp}, \beta)\) is defined as the quotient of an open subset \(Z \subset \A^r\) by the action of a group \(G\) that acts on \(\A^{r}\) through a homomorphism \(G \to \Gm^r\).
  In our case, \(N\) is a free \(\Z\)-module.
  From the construction of \(G \to \Gm^r\) in \cite[\S~2]{bor.che.smi:05}, it follows that \(G \to \Gm^r\) is injective.
  Let \(\mathscr X = [Z/G]\) and \(\overline{\mathscr X} = [Z/\Gm^r]\).
  It is easy to see that we have the pull-back diagram
  \[
    \begin{tikzcd}
      \mathscr X \ar{r}\ar{d} & \overline{\mathscr X}\ar{d} \\
      \A^n\ar{r} & {[\A^{n}/\Gm^n]}.
    \end{tikzcd}
  \]
\end{remark}

Given \(\Lambda\), we use two natural choices of \(\beta\).
For the first, denoted by \(\beta^{\rm can}\), we let \(\beta_F\) be the shortest vector with integer coordinates on the ray \(F^{\perp}\).
For the second, denoted by \(\beta^{\rm res}\), we let \(\beta_F\) be the shortest vector with integer coordinates on the ray \(F^{\perp}\) such that the value of \(\beta_F\) on \(F\) is an integer.
Then we have a map
\[ \mathscr X_{\Lambda, \beta^{\rm res}} \to \mathscr X_{\Lambda, \beta^{\rm can}},\]
which is a sequence of root stacks along the divisors defined by the rays.
Precisely, it is the root stack of order \(\beta^{\rm res}_F / \beta^{\rm can}_F\) along the divisor defined by the ray \(F^{\perp}\).
The map \(\mathscr X_{\Lambda, \beta^{\rm can}} \to X_{\Lambda}\) is called the \emph{canonical desingularisation}.
The map \(\mathscr X_{\Lambda, \beta^{\rm can}} \to \A^n\) is an isomorphism away from the origin.
The map \(\mathscr X_{\Lambda, \beta^{\rm res}} \to \A^n\) is an isomorphism away from the union of the coordinate hyperplanes.

Let \(x_1, \dots, x_n\) be variables.
For \(p = (p_1, \dots, p_n) \in \mathbf{Z}^n_{\geq 0}\), we write \(x^p\) for the monomial \(x_1^{p_1}\cdots x_n^{p_n}\).
A \emph{weighted monomial} is a pair \((x^p,d)\), where \(d\) is a positive integer.
Let \(L\) be a set of weighted monomials.
Let \(\Lambda\) be the rational Newton polyhedron defined by the points \(\frac{1}{d} p\) for \((p,d) \in L\), that is, the convex hull of the union of \(\frac{1}{d} p + \R^n_{\geq 0}\) for \((x^p,d) \in L\).
Assume that \(\Lambda\) is simplicial.
Then we have the stacky blow-ups \(\mathscr X_{\Lambda, \beta^{\rm res}}\) and \(\mathscr X_{\Lambda, \beta^{\rm can}}\).
We call \(\mathscr X_{\Lambda, \beta^{\rm can}}\) the \emph{canonical weighted blow-up} in the set of weighted monomials \(L\), and \(\mathscr X_{\Lambda, \beta^{\rm res}}\) the \emph{resolving weighted blow-up}.
The following two propositions justify the name.
\begin{proposition}\label{prop:monomial}
  In the setup above, let \(v = \frac{1}{d} p\) be a vertex of \(\Lambda\).
  Consider the chart
  \[ \spec \k[u_1, \dots, u_n]  \to \mathscr X_{\Lambda, \beta^{\rm res}}\]
  defined by \(v\).
  The image of \(x^p\) in \(\k[u_1,\dots,u_n]\) is the \(d\)-th power of a monomial \(u\).
  Furthermore, for every \((x^q,e) \in L\), the monomial \(u^e\) divides the image of \(x^q\).
\end{proposition}
\begin{proof}
  Set \(\beta_i = \beta^{\rm res}_{F_i}\).
  Using \eqref{eqn:blowupmap}, we see that
  \begin{align*}
    x^p &\mapsto \prod_i \prod_j u_j^{p_i\langle \beta_j, e_i \rangle} =  \prod_j u_j^{\langle \beta_j, p \rangle}.
  \end{align*}
  By the choice of \(\beta_j\), the quantity \(\langle  \beta_j, p/d \rangle\) is a non-negative integer.
  Thus, \(x^p\) maps to the \(d\)-th power of the monomial
  \[ u = \prod_j u_{j}^{\langle  \beta_j, p/d \rangle}.\]
  Consider \((x^q,e) \in L\).
  Let \(r_1, \dots, r_n\) be the rays of \(\Lambda\) incident to \(v\).
  Then the point \(q/e\) is in the cone defined by the vertex \(v\) and the rays spanned by \(r_1,\dots, r_n\).
  That is, we can write
  \[ q/e = p/d + a_1 r_1 + \cdots + a_nr_n\]
  for some non-negative rational numbers \(a_i\).
  By applying \(\beta_i\) to both sides, we see that \(e \cdot a_i\) is a non-negative integer.
  Using \eqref{eqn:blowupmap} again, we get 
  \[ x^q \mapsto u^e \prod u_i^{e \cdot a_i}.\]
\end{proof}
Let \(d_1, \dots, d_m\) be positive integers and let \(\mathscr P(d_1,\dots,d_m)\) be the weighted projective stack
\[ \mathscr P(d_1, \dots, d_m) = [(\A^m - 0) \big / \Gm],\]
where \(\Gm\) acts coordinate-wise by weights \(d_1, \dots, d_m\).
Consider the rational map \[\A^n \dashrightarrow \mathscr P(d_1, \dots, d_m)\] defined by the monomials \(x^{p_1}, \dots, x^{p_m}\); that is,
\begin{equation}\label{eqn:ratmon}
  (x_1, \dots, x_n) \mapsto [x^{p_1} : \cdots : x^{p_m}].
\end{equation}
Let \(\Lambda\) be the Newton polyhedron defined by the weighted monomials \((x^{p_1},d_1), \dots, (x^{p_m},d_m)\).
Assume that \(\Lambda\) is simplicial.
\begin{proposition}\label{prop:resolution}
  The map \eqref{eqn:ratmon} extends uniquely to a morphism
  \[ \mathscr X_{\Lambda, \beta^{\rm res}} \to \mathscr P(d_1, \dots, d_m).\]
\end{proposition}
\begin{proof}
  The domain is normal and the co-domain is separated.
  Therefore, if the map extends, it extends uniquely \cite[Appendix~A]{fan.man.nir:10}.
  To see that it extends, we may work locally on charts.
  Let \(v\) be a vertex of \(\Lambda\), say \(v = p_i/d_i\).
By \Cref{prop:monomial}, on the chart \(\spec \k[u_1,\dots,u_n]\), the pull-back of \(x^{p_i}\) is \(u^{d_i}\) for a monomial \(u\), and the pull-back of \(x^{p_j}\) is divisible by \(u^{d_j}\).
The extension of \eqref{eqn:ratmon} on \(\spec \k[u_1,\dots, u_n]\) is given by
\[[x^{p_1}u^{-d_{1}}: \cdots : x^{p_{i-1}}u^{-d_{i-1}}:1:x^{p_{i+1}}u^{-d_{i+1}} : \cdots :x^{p_m}u^{-d_m}]. \]
\end{proof}

We reformulate \eqref{prop:resolution} to suit our setting.
\begin{corollary}\label{cor:resolution}
  Let \(W_1, \dots, W_{m}\) be finite dimensional \(\k\)-vector spaces and set \(W = \bigoplus_i W_i\).
  Let \(\sP W\) be the weighted projective stack where \(W_i\) has weight \(d_i > 0\).
  Let \(f \colon \A^n \dashrightarrow \sP W\) be the rational map defined by \(f_i \in W_i \otimes \A[x_{1}, \dots, x_n]\) and assume that the coordinates of \(f_{i}\) generate the monomial ideal \(\langle x^{p_i} \rangle\).
  Let \(\Lambda\) be the Newton polyhedron defined by the weighted monomials \((x^{p_1},d_1), \dots, (x^{p_m},d_m)\).
  Then the rational map \(f\) extends uniquely to a morphism
  \[ \mathscr X_{\Lambda, \beta^{\rm res}} \to \sP W.\]
\end{corollary}
\begin{proof}
  We follow the proof of \Cref{prop:resolution}.
  Let \(v\) be a vertex of \(\Lambda\), say \(v = p_i/d_i\).
  Let \(u^{d_i}\) be the pull-back of \(x^{p_i}\) to the chart defined by \(v\).
  Then for all \(j\), the element \(u^{-d_j}f_j \in W_j \otimes \k[u_1,\dots,u_m]\) has polynomial coordinates.
  Furthermore, for \(j = i\), the coordinates generate the unit ideal.
  On this chart, the extension of the rational map \(f\) is given by \([u^{-d_1}f_1: \cdots : u^{-d_m}f_m]\).
\end{proof}

\section{Class of an orbit using a complete parametrisation}\label{sec:param}
Let \(W\) be a finite dimensional representation of \(\GL({m})\).
Consider the central \(\G_{m} \subset \GL({m})\), and assume that it acts on \(W\) by positive weights.
We denote by \(\sP W\) the weighted projective stack
\[ \sP W = [W - 0 \big / \G_{m}].\]

Fix a non-zero vector \(w \in W\), and let
\[ [w] \colon \spec \k \to \sP W\]
be the corresponding point of \(\sP W\).
By the \emph{stabiliser} \(\Gamma\) of \(w\), we mean the fiber product
\begin{equation}\label{eq:stab}
\begin{tikzcd}
  \Gamma \ar{r}\ar{d} & \PGL(m)  \ar{d}{a \mapsto a\cdot w} \\
  \spec \k \ar{r}{[w]} & \sP W.
\end{tikzcd}
\end{equation}
We have the diagram
\[
  \begin{tikzcd}
     \Gamma\ar{d}\ar{r}  &\GL(m) \ar{r}\ar{d}{a \mapsto a\cdot w}&\PGL_m\ar{d}{a \mapsto a \cdot w}\\
    \spec \k \ar{r} & W \ar{r} & \sP W 
  \end{tikzcd}
\]
in which the right square and the outer square are cartesian.
Therefore, the left square is also cartesian.
Therefore, \(\Gamma\) is simply the stabiliser of \(w\) in \(\GL(m)\).

A \emph{complete orbit parametrisation} of \([w]\) is a proper morphism
\[ i \colon X \to \sP W,\]
where \(X\) is a Deligne--Mumford stack together with the action of \(\PGL(m)\) and \(i\) is a \(\PGL(m)\)-equivariant map such that there exists an open subscheme \(U \subset X\) isomorphic to \(\PGL(m)\) as a \(\PGL(m)\)-scheme and a point \(x \in U\) whose image is \([w]\).
The \emph{orbit of \([w]\)}, denoted by \(\Orb([w])\), is the Zariski closure in \(\sP W\) of \(\PGL(m) \cdot [w]\), with the reduced scheme structure.
\begin{proposition}\label{prop:weighted-complete-param}
  Let \(i \colon X \to \sP W\) be a complete parametrisation of the orbit of \(w\).
  Assume that the stabiliser \(\Gamma \subset \GL(m)\) of \(w\) is finite.
  Then, we have the equality of cycles
  \[ i_{*}[X] = |\Gamma| [\Orb([w])].\]
\end{proposition}

\begin{proof}
  We have the fiber product
\begin{equation}\label{eq:staborb}
    \begin{tikzcd}
      \Gamma \ar{r}\ar{d} & \PGL(m)\ar{d} \\
      \spec \k \ar{r}{[w]} & \Orb([w]).
    \end{tikzcd}
  \end{equation}
  Consider the open inclusion \(\PGL(m) \to X\) that sends \(a\) to \(a \cdot x\).
  The image of this inclusion is \(U\).
  The points of \(X\) in the complement of \(U\) are stabilised by a positive dimensional subgroup of \(\PGL(m)\) and hence they map to points in \(\Orb([w])\) that are stabilised by a positive dimensional subgroup.
  In particular, they do not map to \([w]\).
  As a result, the fiber product \eqref{eq:staborb} gives the fiber product
    \[
    \begin{tikzcd}
      \Gamma \ar{r}\ar{d} & X \ar{d} \\
      \spec \k \ar{r}{[w]} & \Orb([w]).
    \end{tikzcd}
  \]
  We see that the map \(X \to \Orb([w])\) is generically finite of degree \(|\Gamma|\).
  The proposition follows.
\end{proof}

We now give a cohomological formula for the push-forward.
We first need a lemma, adapted from \cite[Proposition~2.1]{deo.pat:22}.
Let \(U\) be a vector space of dimension \(N\) with the action of an algebraic group \(G\).
Set \(U^{*} = U-0\) and let \(\pi \colon U^{*} \to \P U\) be the projection.
\begin{lemma}\label{prop:pushforward}
  Let \(Y\) be a Deligne--Mumford stack with a \(G\)-action and a \(G\)-equivariant map \(\phi \colon Y \to \P U\).
  Then, in \(A_G(U^{*})\), we have the equality
  \[ \pi^{*} \phi_{*}[Y] =  \int_{Y} \frac{c_{N}(U)}{\phi^* c_1 \O(-1)}.\]
\end{lemma}
  The integral on the right is the push-forward \(A_G(Y) \to A_G\),  considered as an element of \(A_G(U^{*})\) via the pull-back \(A_G \to A_G(U^{*})\).
\begin{proof}
  Let \(Q\) be the cokernel of \(\phi^{*}\O(-1) \to U \otimes \O_Y\).
  On \(Y \times \P U\), let \(\pi_{i}\) for \(i = 1,2\) be the two projections.
  The vanishing locus of the composite map
  \[ \pi_2^{*} \O(-1) \to U \otimes \O_{Y \times \P U} \to \pi_1^{*} Q\]
  is precisely the graph \(Z\) of \(\phi \colon Y \to \P U\).
  Therefore, we have
  \begin{equation}\label{eqn:gamma}
    [Z] = c_{N-1}(\pi_1^{*}Q \otimes \pi_2^{*}\O(1)) [Y \times \P U].
  \end{equation}
 Consider the fiber square
  \[
    \begin{tikzcd}
      X \times U^{*} \ar{r}{\widetilde\pi_2}\ar{d}{\widetilde\pi}& U^{*}\ar{d}{\pi} \\
      X \times \P U \ar{r}{\pi_2} & \P U.
    \end{tikzcd}
  \]
  By the push-pull formula, we have
  \begin{equation}\label{eqn:pushpull}
    {\widetilde{\pi_2}}_{*}\widetilde\pi^{*} [Z] = \pi^{*} {\pi_2}_{*}[Z].
  \end{equation}
  The right-hand side of \eqref{eqn:pushpull} is \(\pi^{*}\phi_{*}[X]\).
  Since the pull-back of \(\O_{\P U}(1)\) to \(U^{*}\) is trivial, \eqref{eqn:gamma} shows that 
  \[\widetilde\pi^{*}[Z] = c_{N-1}(\pi_1^{*}Q)[Y \times U^{*}].\]
  The statement follows by applying \({\widetilde{\pi_2}}_{*}\) to the above equation.
\end{proof}

We need an analogue of \Cref{prop:pushforward} for weighted projective spaces.
Let \(W\) be a vector space of dimension \(N\) with an action of a torus \(T\).
Set \(W^{*} = W - 0\) and \(\sP W = [W^{*} \big / \G_m]\), where \(\G_m\) acts on \(W\) by positive weights and this action commutes with the action of \(T\).
Let \(\pi \colon W^{*} \to \sP W\) be the projection.
\begin{lemma}\label{prop:weightedpushforward}
  Let \(X\) be a Deligne--Mumford stack with a \(T\)-action and a \(T\)-equivariant map \(\phi \colon X \to \sP W\).
  Then, in \(A_T(W^{*})\), we have the equality
  \[ \pi^{*} \phi_{*}[X] = \int_X \frac{c_N(W)}{\phi^* c_1 \O_{\sP W}(-1)}.\]
\end{lemma}
We understand the right-hand side in the same sense as in \Cref{prop:pushforward}.
\begin{proof}
  It suffices to prove the equality in \(A_{\widetilde T}(W^{*})\) where \(\widetilde T \to T\) is a finite cover by another torus.
  Choose a basis \(\langle  w_i \rangle\) of \(W\) compatible with the action of \(T\) and \(\Gm\).
  Suppose \(T\) acts on \(w_i\) by the character \(\chi_i \in \Hom(T, \G_m)\) and \(\G_m\) acts on \(w_i\) by weight \(d_i\).
  Let \(\widetilde T \to T\) be a finite cover by a torus such that the image of \(\chi_{i}\) in \(\Hom(\widetilde T, \Gm)\) is divisible by \(d_i\).
  Let \(U\) be the \(\k\)-span of the symbols \(u_{i}\).
  Equip \(U\) with a \(\widetilde T\) action so that \(\widetilde T\) acts on \(u_{i}\) by the character \(\frac{1}{d_i} \chi_{i}\) and with the \(\Gm\) action by weight \(1\).
  Then the map
  \[ \mu \colon U \to W\]
  defined by \(\sum x_iu_i \mapsto \sum x_i^{d_i} w_i\) is equivariant for the \(\widetilde T\) and \(\G_{m}\) actions and finite of degree
  \[ \deg \mu = \prod d_i.\]
  Under the induced map \(\mu \colon \P U \to \sP W\), the pull-back of \(\O_{\sP W}(-1)\) is \(\O_{\P U}(-1)\).
  Define \(Y\) by the pull-back diagram
  \[
  \begin{tikzcd}
    Y \ar{r}{\widetilde\phi}\ar{d}{\widetilde\mu} & \P U \ar{d}{\mu} \\
    X \ar{r}{\phi} &\sP W.
  \end{tikzcd}
\]
Set \(U^{*} = U-0\) and denote by \(\widetilde\pi \colon U^{*} \to \P U\) the projection.
By \Cref{prop:pushforward}, in \(A_{\widetilde T}(U^{*})\) we have
\begin{equation}\label{eqn:topintegral}
  \widetilde\pi^{*} \widetilde\phi_{*}[Y] = \int_Y \frac{c_N(U)}{\phi^* c_1 \O_{\P U}(-1)}.
\end{equation}
Note that \(c_{N}(U) = \prod d_i^{-1} c_N(W)\).
Since \(Y \to X\) is of degree \(\prod d_i\), the integral on the right-hand side of \eqref{eqn:topintegral} is equal to
\[\int_X \frac{c_N(W)}{\phi^{*}c_1 \O_{\sP W}(-1)}.\]
Now the statement follows by applying \(\mu_{*} \colon A_{\widetilde T}(U^{*}) \to A_{\widetilde T}(W^{*})\) to both sides of \eqref{eqn:topintegral}.
\end{proof}

Let \(\Orb(w) \subset W\) be the closure of the \(\GL(m)\)-orbit of \(w\), with the reduced scheme structure.
Let \(N = \dim W\).
\begin{proposition}
  \label{prop:weightedformula}
  Let \(i \colon X \to \sP W\) be a complete parametrisation of the orbit of \(w \in W\).
  Assume that the stabiliser \(\Gamma \subset \GL(m)\) of \(w\) is finite.
  Then, in \(A_{G}(W) = A_{G}\), we have
  \[
    |\Gamma| \cdot [\Orb(w)] = \int_{X} \frac{c_{N}(W)}{i^{*}c_1\O_{\sP W}(-1)}.
  \]
\end{proposition}
The integral on the right is the push-forward \(A_G(X) \to A_G\).
\begin{proof}
  Let \(T \subset \GL(m)\) be a maximal torus.
  It suffices to prove the equality in \(A_T(W)\).
  \Cref{prop:weighted-complete-param} and \Cref{prop:weightedpushforward} together give the equality in \(A_T(W^{*})\).
  But then the equality also holds in \(A_T(W)\) since \(A^i_T(U) = A^i_T(U^{*})\) for \(i = \codim \Orb(w)\).
\end{proof}

\section{Twist invariance}\label{sec:twist}
Consider a representation \(W\) of \(\GL(2)\) defined by
\[ \rho \colon \GL(2) \to \GL(W).\]
For \(n \in \Z\), we have the surjective homomorphism
\[ T_n \colon \GL(2) \to \GL(2), \quad M \mapsto M \cdot (\det M)^{n},\]
whose kernel is the diagonally embedded \(\mu_{n+2} \subset \GL(2)\).
The composite \(\rho \circ T_n\) gives a new representation \(\GL(2) \to \GL(W)\) which we call \(W(n)\).
Note that 
\begin{equation}\label{eqn:twist}
 \text{if } W \cong \Sym^{a-b}V \otimes \det V^b, \text{ then } W(n) \cong \Sym^{a-b} V \otimes \det{}^{b + n(a+b)} V.
\end{equation}
Observe that the identity map \(W(n) \to W\) together with \(T_n \colon \GL(2) \to \GL(2)\) induces a map
\[ e_n \colon [W(n)/\GL(2)] \to [W/\GL(2)].\]

Given \(w \in W\), the \(\GL(2)\)-orbit closure of \(w\) in \(W\) under \(\rho\) is equal to that of \(w\) in \(W(n)\).
But to distinguish the ambient representations, we denote them by \(\Orb(w)\) and \(\Orb(w)(n)\), respectively.
Then
\[ \Orb(w)(n) = e_n^{-1} (\Orb(w)),\]
and hence
\[ [\Orb(w)(n)] = e_n^{*} ([\Orb(w)]) \in A_{\GL(2)}.\]
The map
\[ e_n^{*} \colon A_{\GL(2)} \to A_{\GL(2)}\]
is easy to describe.
Thinking of \(A_{\GL(2)}\) as the subring of \(\Q[v_1,v_2]\) consisting of symmetric polynomials, it is given by
\begin{equation}\label{eqn:twistpull}
  e_n^{*} \colon v_1 \mapsto v_1 + n(v_1+v_2) \text{ and } v_2 \mapsto v_2 + n(v_1+v_{2}).
\end{equation}

Let \(\Gamma\) and \(\Gamma(n)\) be the stabilisers of \(w\) under \(\rho\) and \(\rho \circ T_n\), respectively.
Then we have the sequence
\[ 1 \to \mu_{n+2} \to \Gamma(n) \to \Gamma \to 1.\]
In particular, we have \(|\Gamma(n)| = (n+2)|\Gamma|\).

Given \(u \in \P^1\), let \(\Lambda^u\) be the Newton polygon associated to \(w \in W\) and \(\Lambda^u(n)\) the Newton polygon associated to \(w \in W(n)\).
Using \eqref{eqn:twist}, it follows that \(\Lambda^u(n) \subset \R^2\) is obtained from \(\Lambda^u \subset \R^2\) by applying the transformation
\begin{equation}\label{eqn:shear}
  (x,y) \mapsto \frac{1}{n+2}(x+1,y+1).
\end{equation}
Let \(Q\) be the polynomial on the right-hand side of \eqref{eqn:main} in the main theorem for \(w \in W\) and \(Q(n)\) the corresponding polynomial for \(w \in W(n)\).
Using \eqref{eqn:twistpull} and \eqref{eqn:shear}, it is easy to check that
\[ e_n^{*}(Q) = Q(n)/(n+2).\]
Since we also have \[e_n^{*}(|\Gamma|[\Orb(w)]) = |\Gamma(n)|[\Orb(w)(n)]/(n+2),\] the main theorem holds for \(w \in W\)  if and only if it holds for \(w \in W(n)\).
Thus, in the proof, we are free to replace \(W\) by \(W(n)\) for any \(n\).
In particular, by choosing a sufficiently large \(n\), we may assume without loss of generality that 
\begin{equation}\label{eq:positiveb}
  W \cong \bigoplus \Sym^{a_i-b_i}V \otimes \det{}^{b_i} V \text{ with } a_i \geq b_i \geq 0.
\end{equation}

\section{Complete orbit parametrisations of \(\GL(2)\)-orbits}\label{sec:completeparam}
Recall that we have a 2-dimensional vector space \(V\) and 
\[ W = W_{1} \oplus \cdots \oplus W_{n}, \text{ where } W_i = \Sym^{a_i-b_i}V \otimes \det{V}^{b_i}.\]
We set \(d_i = a_i+b_i\), which we call the \emph{weight} of \(W_i\), and assume \(d_i > 0\) for all \(i\).
Also assume that \(b_i \geq 0\); this can be achieved after twisting \(W\) as in \Cref{sec:twist}.
Consider the central \(\mathbf{G}_{m} \to \GL V\) given by \(t \mapsto t \cdot I\).
Observe that \(t \in \G_{m}\) scales the elements of \(W_i\) by \(t^{d_i}\).
If \(U\) is another 2-dimensional vector space, then by \(W_{i}(U)\) we mean the representation
\[ W_{i}(U) = \Sym^{a_{i}-b_{i}}U \otimes \det{}^{b_{i}}U,\]
and by \(W(U)\) the direct sum
\[ W(U) = \bigoplus_{i} W_{i}(U).\]
Let \(\sP W\) be the weighted projective stack
\[ \sP W = \left[ W - 0 \big / \G_{m} \right]\]
Let \(U\) be another two-dimensional vector space and set
\[ M = \P \Hom(U, V).\]
Fix a non-zero \(w \in W(U)\).
Let \(w_{i} \in W_{i}(U)\) be the \(i\)-th component of \(w\).

Let \(I = \{i \mid w_i \neq 0\}\) and \(J = \{1,\dots, n\} - I\).
Set \(W_I = \oplus_{i \in I} W_i\) and similarly for \(W_J\).
Let \(w_{I}\) be the projection of \(w\) to \(W_I\).
Plainly, we have \(\Orb(w) = \Orb(w_I) \times \{0\} \subset W_I \oplus W_J\), and hence
\[[\Orb(w)] = c_{\dim W_J}(W_J)[\Orb(w_I)].\]
Using this, we see that it suffices to prove the main theorem when \(J = \emptyset\).
So, assume that \(w_i \neq 0\) for all \(i\).

We have a rational map
\begin{equation}\label{eq:ratmap}
  M \dashrightarrow \sP W
\end{equation}
defined by
\[ m \mapsto [ m \cdot w].\]
It is defined on the locus of \(m\) such that \(m \cdot w \neq 0\).
More formally, on \(M = \P\Hom(U, V)\), we have the universal homomorphism
\[
  e \colon U \otimes \O_M(-1) \to  V \otimes \O_M,
\]
which induces
\[ W_{i}(U) \otimes \O_M(-d_{i}) \to W_{i}(V) \otimes \O_M.\]
By pre-composing with the section
\[ w_i \colon \O_M \to W_i(U_i) \otimes \O_M,\]
we get the map
\begin{equation}\label{eq:irrmap}
  \O_M(-d_{i}) \to W_{i}(V) \otimes \O_M.
\end{equation}
The maps in \eqref{eq:irrmap} define a morphism to \(\sP W\) on the open subset of \(M\) where at least one of the maps is non-zero.
Observe that this open subset includes all points of \(M\) corresponding to invertible homomorphisms.

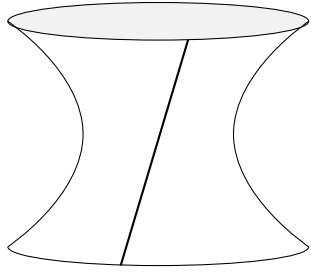
\begin{figure}
    \centering
  \begin{tikzpicture}
    \draw[thick]
    (-0.5, -1.75) -- (0.5,1.6);
    \draw[fill=black!5!white] (0,1.5) ellipse (2cm and 0.25cm);
    \draw (2,-1.5) arc (0:-180:2cm and 0.25cm);
    \draw plot [smooth, tension=1] coordinates {(-2,1.5) (-1,0) (-2,-1.5)};
    \draw plot [smooth, tension=1] coordinates {(2,1.5) (1,0) (2,-1.5)};
  \end{tikzpicture}
  \caption{The scheme theoretic zero locus of the map \eqref{eq:irrmap} is cut out locally by an ideal of the form \(I_{K_u}^r \cdot I_{\Delta}^{b}\), where \(\Delta \subset M\) is the determinant quadric and \(K_u \subset \Delta\) are certain lines on it.}
  \label{fig:kd}
\end{figure}

We now describe the scheme theoretic zero locus of the map \eqref{eq:irrmap}.
It is supported on the determinant quadric
\[ \Delta = \{m \in M \mid \det m = 0\},\]
and it has embedded primes supported on lines of one ruling of this quadric (see \Cref{prop:zerolocus} and \Cref{fig:kd}).
Given a point \(u \in \P U\), let \(K_{u} \subset M\) be the line defined by
\[ K_{u} = \{m \in M \mid m u = 0\}.\]
Observe that as \(u\) varies in \(\P U \cong \P^1\), the lines \(K_u\) sweep out one of the two rulings of \(\Delta\).

\begin{proposition}\label{prop:zerolocus}
  Suppose \(w_i = f \otimes \delta^{b_i}\), where \(f \in \Sym^{a_i-b_i}(U)\) and \(\delta \in \det U\) are non-zero.
  Take \(m \in \Delta \subset M\)  and let \(u \in \P U\) be the kernel of \(m\).
    Suppose \(f\) vanishes to order \(r\) at \(u\).
  Then, in a neighbourhood of \(m\), the scheme theoretic zero locus of the map 
  \[ e \colon \O_M(-d_i) \to W_i(V) \otimes \O_M,\]
  defined in \eqref{eq:irrmap}, is cut out by the ideal \(I_{K_{u}}^r \cdot I_{\Delta}^{b_i}\).
\end{proposition}
\begin{proof}
 Since \(i\) is fixed, we omit it from the subscript in \(a_i, b_i,\) and \(d_i\), and do a local calculation.
  Denote by \(u_{2} \in U\) a lift of \(u \in \P U\).
  We choose a linearly independent vector \(u_{1} \in U\) and take \((u_{1},u_{2})\) as a basis of \(U\).
  In this basis, the point \(u \in \P U\) is given by \([0:1]\).

  Choose a basis \((v_{1},v_{2})\) of \(V\) and suppose that in the chosen bases, the map \(m \colon U \to V\)  is given by
  \[ m = \begin{pmatrix} 1 & 0 \\ 0 & 0 \end{pmatrix}.\]
   We can take \(\delta = u_1 \wedge u_2\).
  Consider the affine neighbourhood of \(m \in M\) given by matrices of the form
  \[ \begin{pmatrix} 1 & x \\ z & y+xz \end{pmatrix}.\]
  Up to multiplication by a non-zero scalar, the element \(f = f(u_1,u_2)\) has the form
  \begin{equation}\label{eq:p}
    f(u_1,u_2) = u_{1}^{a-b-r}u_{2}^{r} + * \cdot u_{1}^{a-b-r-1}u_{2}^{r+1} + \cdots.
  \end{equation}
  Substituting \(u_{1} \mapsto v_{1} + z v_{2}\) and \(u_{2} \mapsto xv_{1}+(y+xz)v_{2}\) yields
  \begin{equation}\label{eqn:pex}
    \begin{split}
      (Mf)(v_1,v_{2}) &= (v_1+zv_2)^{a-b-r}(xv_{1}+(y+xz)v_{2})^{r}  \\
      & \qquad + * \cdot (v_{1}+zv_{2})^{a-b-r-1}(xv_{1}+(y+xz)v_{2})^{r+1}+\cdots,
    \end{split}
  \end{equation}
  and,
  \begin{equation}\label{eqn:pexdet}
    (M \delta)(v_1,v_2) = y (v_1 \wedge v_2).
  \end{equation}
  Observe that the ideal generated by the coefficients of \(Mf\) is \(\langle  x,y \rangle^{r}\).
  The ideal \(\langle  x,y \rangle\) is precisely the ideal \(I_{K_{u}}\) and the ideal \(\langle y \rangle\) is precisely the ideal \(I_{\Delta}\).
  So the ideal generated by the coefficients of \(M(f \otimes \delta^b)\) is \(I_{K_u}^{r} \cdot I_{\Delta}^{b}\), as required.
\end{proof}

We now resolve the rational map \(M \dashrightarrow \sP W\) using a stacky blow-up.
Let the components of \(w\) be \(w_i = f_i \otimes \delta^{b_i}\).
Let \(A \subset \P U\) be any finite set that includes the common zeros of \(f_i\) for \(i\) that realise the minimum \(\min_i b_i/d_i\).
Let \(\Bl_A M \to M\) be the blow-up of \(M\) along the lines \(K_u\) for \(u \in A\).
Let \(E_u \subset \Bl_A M\) be the exceptional divisor over \(K_u\) and \(D \subset \Bl_{A}M\) the proper transform of \(\Delta \subset M\).
Let \(r_i^u\) be the order of vanishing of \(f_i\) at \(u\).
Let 
\[ e \colon \O_{\Bl_A M}(-d_i) \to W_i(V) \otimes \O_{\Bl_A M}\]
be the pull-back of \eqref{eq:irrmap}.
By \Cref{prop:zerolocus}, the ideal generated by the components of this map is \(I_{E_u}^{r_i^u+b_i}I_D^{b_i}\).

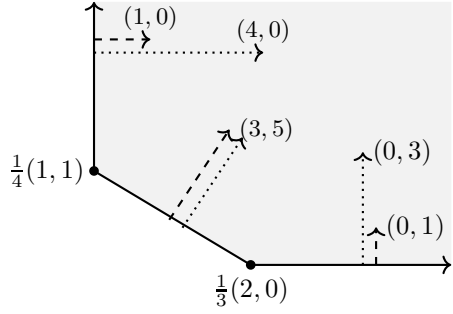
\begin{figure}
  \begin{tikzpicture}[scale=5, thick]
    \draw[draw=none,fill=black!5!white] (1/4,.7) -- (1/4,1/4) -- (2/3,0) -- (1.2,0) -- (1.2,.7) -- (1/4,.7);
    \draw[->] (2/3,0) -- (1.2,0);
    \draw[->] (1/4,1/4) -- (1/4,.7);
    \draw (2/3,0) -- (1/4,1/4);
    \draw[fill]
    (2/3, 0) circle (0.01) node[below] {\(\frac{1}{3}(2,0)\)}
    (1/4,1/4) circle (0.01) node[left] {\(\frac{1}{4}(1,1)\)};
    \draw[dashed, ->] (1/4,.6) -- (.4,.6) node[above] {\small \((1,0)\)};
    \begin{scope}[yshift=-1]
      \draw[dotted, ->] (1/4,.6) -- (.7,.6) node[above] {\small \((4,0)\)};
    \end{scope}
    \draw[dashed, ->] (1,0) -- (1,.1) node[right] {\((0,1)\)};
    \begin{scope}[xshift=-1]
      \draw[dotted, ->] (1,0) -- (1,.3) node[right] {\((0,3)\)};
    \end{scope}
    \draw[dashed, ->] (.45,.12) -- (.61,.36) node[right] {\small \((3,5)\)};
    \begin{scope}[xshift=1, yshift=-3/5]
      \draw[dotted, ->] (.45,.12) -- (.61,.36);
    \end{scope}
  \end{tikzpicture}
  \caption{
    For \(W = \Sym^3V \oplus \Sym^2 V \otimes \det V\), the Newton polygon \(\Lambda^u\) at \(u \in \P^1\) with the vanishing orders \(r_1 = 2\) and \(r_2 = 0\).
    The short normal vectors (dashed) represent \(\beta^{\rm can}\) and the longer ones (dotted) represent \(\beta^{\rm res}\).
  }
  \label{fig:blowup}
\end{figure}

Fix \(u \in A \subset \P U\).
Let \(\Lambda = \Lambda^u \subset \R^2\) be the Newton polygon defined by the set of weighted monomials \(\{(x^{r_i^u + b_i}y^{b_i}, d_i) \mid i = 1, \dots, n\}\).
Let \(\overline {\mathscr X}_{\Lambda, \beta} \to [\A^2/\Gm^2]\) be the blow-up defined in \Cref{sec:stackyblowup} (see \Cref{rem:eqvblowup}) for \(\beta = \beta^{\rm can}\)  and \(\beta = \beta^{\rm res}\).
Let \(\Bl_AM \to [\A^2/\Gm^2]\) be the map defined by the divisors \(E_u\) and \(D\) and let
\(\mathscr M^{\rm res}_{u} \to \Bl_AM\) and \(\mathscr M^{\rm can}_u \to \Bl_AM\) be the pullbacks of \(\overline {\mathscr X}_{\Lambda, \beta^{\rm res}} \to [\A^2/\Gm^2]\) and \(\overline {\mathscr X}_{\Lambda, \beta^{\rm can}} \to [\A^2/\Gm^2]\).
Since \(E_u\) and \(D\) are smooth, normal crossings divisors, the map \(\Bl_AM \to [\A^2/\Gm^2]\) is smooth, and therefore both \(\mathscr M^{\rm res}_u\) and \(\mathscr M^{\rm can}_u\) are smooth.
See \Cref{fig:blowup} for an example of \(\Lambda^{u}\) with \(\beta^{\rm can}\) and \(\beta^{\rm res}\).

We can write local charts for \(\mathscr M^{\rm res}_u\) and \(\mathscr M^{\rm can}_u\) by simply substituting local equations of \(E_u\) and \(D\) in the local charts described in \Cref{sec:stackyblowup}.
Let \(\lambda = (\lambda_1,\lambda_2)\) be the vertex of the Newton polyhedron \(\Lambda = \Lambda^u\) with the smallest second coordinate.
Then \(\lambda_2 = \min_i(b_i/d_i)\).
Suppose \(\lambda_2 = b/d\) where \(\gcd(b,d) = 1\).
Note that \(\lambda + \R_{\geq 0} \times 0\) is a ray of \(\Lambda\).
Its associated divisor is the vanishing locus of \(y\), which pulls back to \(D \subset \Bl_A M\).
The functional \(\beta_1^{\rm can}\) associated to this ray is the projection \(p_2\colon (a,b) \mapsto b\).
On the other hand, the functional \(\beta_1^{\rm res}\) is \(d \cdot p_2\).
As a result, \(\mathscr M^{\rm res}_u \to \Bl_AM\) over the complement of \(E_u\) is the root stack along \(D\) of order \(d\).
Note that \(d\) is independent of \(u \in A\).

Let \(\mathscr M^{\rm res} \to \Bl_AM\) and \(\mathscr M^{\rm can} \to \Bl_AM\) be the blow-ups as above carried out for all \(u \in A\) at once.
That is, for all \(u \in A\), in a neighbourhood of \(E_u\), the map \(\mathscr M^{\rm res} \to \Bl_AM\) is the blow-up \(\mathscr M^{\rm res}_u \to \Bl_AM \), and similarly for \(\mathscr M^{\rm can} \to \Bl_AM\).
We have maps
\[ \mathscr M^{\rm res} \to \mathscr M^{\rm can} \to \Bl_AM.\]
The map \(\mathscr M^{\rm can} \to \Bl_A M\) is an isomorphism away from the union of the lines \(E_u \cap D\) for \(u \in A\).
The map \(\mathscr M^{\rm res} \to \Bl_A M\) is an isomorphism away from the union of the divisors \(E_u\) for \(u \in A\) and \(D\).
Over the complement of the union of \(E_u\) for \(u \in A\), it is the root stack of order \(d\) along \(D\).

Let \(\mathscr M\) be \(\mathscr M^{\rm can}\) or \(\mathscr M^{\rm res}\).
For every \(g \in \GL V\), it is easy to check that the action map \(g \colon M \to M\) lifts to a morphism \(g \colon \mathscr M \to \mathscr M\).
For \(g,h \in \GL V\), the two morphisms \(h \circ g\) and \(gh\) agree on a dense open subscheme in \(\mathscr M\).
Since \(\mathscr M\) is normal and separated, \cite[Appendix~A]{fan.man.nir:10} implies that there exists a unique \(2\)-morphism \(h \circ g \implies gh\).
As a result, the maps \(g \colon \mathscr M \to \mathscr M\) for \(g \in \GL V\) give an action of \(\GL V\) on \(\mathscr M\).

\begin{proposition}\label{prop:Mresolution}
  The rational map
  \[ M \dashrightarrow \sP W\]
  extends to a regular map
  \[ \iota \colon \mathscr M^{\rm res} \to \sP W,\]
  which is a complete orbit parametrisation of the orbit of \([w] \in \sP W\).
\end{proposition}
\begin{proof}
  The extension exists due to \Cref{prop:resolution} (see \Cref{cor:resolution}).
  It is immediate that \(\iota\) gives a complete orbit parametrisation.
\end{proof}

\section{Atiyah--Bott localisation}\label{sec:localisation}
\Cref{prop:weighted-complete-param} gives a formula for \([\Orb(w)]\) as an integral.
We compute the integral in \Cref{prop:weighted-complete-param} using the Atiyah--Bott localisation formula for stacks \cite[\S~5.3]{kre:99}.
In this section, we use \(\mathscr M\) to denote either \(\mathscr M^{\rm res}\) or \(\mathscr M^{\rm can}\).
A claim about \(\mathscr M\) is understood to hold for both \(\mathscr M^{\rm res}\) and \(\mathscr M^{\rm can}\).
Most such claims will be on the level of points or rational Chow groups, both of which are identical for the two stacks.

Fix a basis \((v_{1}, v_{2})\) of \(V\).
Let \(T \subset \GL V\) be the diagonal torus with respect to the chosen basis.
The \(T\)-fixed locus in \(M\) is the disjoint union of the two lines \(L_i\) for \(i = 1,2\) defined by
\[ L_i = \left\{ m \mid \operatorname{Image}(m) \subset \langle  v_i \rangle \right\}.\]
These are lines on \(\Delta\) of the opposite ruling compared to the lines \(K_u\) (see \Cref{fig:kd}).

\subsection{Fixed points of the \(T\)-action on \(\mathscr M\)}\label{sec:fp}
Let \(\mathscr L^{\rm res}_i \subset \mathscr M^{\rm res}\) and \(\mathscr L^{\rm can}_i \subset \mathscr M^{\rm can}\) be the proper transforms of \(L_i \subset M\) (with the reduced scheme structure).
We use \(\mathscr L_i \subset \mathscr M\) to refer to either one of these.

Fix a \(u \in A \subset \P U\).
Choose a basis \(u_1, u_2\) of \(U\) such that \(u = [u_2]\).
Consider the affine open chart \(\A_{x,y,z}^{3} \subset M\) consisting of matrices of the form
\[
  \begin{pmatrix}
    1 & x \\
    z & y + xz
  \end{pmatrix}.
\]
In this basis, the line \(L_1\), the line \(K_{u}\), and the determinant \(\Delta\) are cut out by
\begin{align*}
  L_{1} &: z = y = 0, \\
  K_{u} &: x = y = 0, \\
  \Delta &: y = 0.
\end{align*}
The line \(L_{2}\) is absent from this chart.
Thinking of \(x\), \(y\), \(z\) as regular functions on this chart, we see that the \(T\)-action is given by
\begin{equation}\label{eqn:t-action-chart}
  (t_1,t_2) \colon (x,y,z) \mapsto (x,t_1t_2^{-1}y, t_1t_2^{-1}z).
\end{equation}
The blow-up \(\Bl_{K_u} M\) has the local description
\[ \{(x,y,z,[X:Y]) \mid Xy=xY\} \subset \spec \k[x,y,z] \times \P^1.\]
On the blow-up, the proper transform of \(L_1\) is cut out by \(z = 0\) and \(Y = 0\).
The only \(T\)-fixed points on the blow-up are the points of the proper transform of \(L_1\) and the point \(p_1^u\) with coordinates \(((0,0,0),[0:1])\).

The proper transform of \(\Delta\) is defined by \(Y = 0\), and is thus contained in the affine chart of the blow-up given by \(X \neq 0\).
This chart is given by
\[ \{(x,y,z,[1:Y]) \mid y = xY\} \cong \spec \k[x,Y,z].\]
The \(T\)-action is given by
\[ (t_1,t_2) \colon (x,Y,z) \mapsto (x, t_{1}^{-1}t_2Y, t_1^{-1}t_2z).\]
The stacky blow-up of this chart is defined by the weighted monomials \((x^{r_i^u+b_i}Y^{b_i}, d_i)\).
Let \(\Lambda^u\) be the Newton polyhedron defined by these weighted monomials (see \Cref{sec:stackyblowup}).
Since \(z\) is absent from the monomials, we may think of \(\Lambda^u\) as a subset of \(\R^2\).
Let \(\lambda(0), \dots, \lambda(k)\) be the vertices of \(\Lambda^u\) arranged from the bottom-right to the top-left.
That is, using subscripts to denote first and second coordinates, we have
\[ \lambda(0)_1 > \cdots > \lambda(k)_1 \text{ and } \lambda(0)_2 < \cdots < \lambda(k)_2.\]
Note that the point corresponding to the bottom-right vertex \(\lambda(0)\) lies on the proper transform of \(L_1\).
It is easy to check that the only \(T\)-fixed points on the stacky blow-up of this chart are:
\begin{enumerate}
\item points of the proper transform of \(L_1\),
\item points corresponding to the vertices \(\lambda(1), \dots, \lambda(k)\).
\end{enumerate}
For \(j = 1,\dots,k\), we label the point corresponding to \(\lambda(j)\) as \(p_{1,j}^u\).

We have analogous points \(p_2^u\) and \(p_{2,j}^u\) over the line \(L_2 \subset M\).

Summarising the discussion above, we see that the \(T\)-fixed locus of \(\mathscr M\) is the disjoint union of 
\begin{enumerate}
\item \(\mathscr L_1 \sqcup \mathscr L_2\)
\item \(\{p_1^u, p_2^u\}\) for \(u \in A\).
\item \(\{p_{1,j}^u, p_{2,j}^u \mid j = 1, \dots, k = k^u\}\) for \(u \in A\).
\end{enumerate}

\subsection{Ingredients of the localisation formula}\label{sec:locingr}
Recall that we have the map \[\iota \colon \mathscr M^{\rm res} \to \sP W,\] which is the complete orbit parametrisation of \([w]\).
We describe the pull-back of \(\O(-1)\) and the normal bundles to the components of the fixed locus of the \(T\)-action as elements of the corresponding (rational) \(T\)-equivariant Grothendieck groups.

Let \(M_T = \Hom(T, \Gm) \otimes \Q\) and \(K_T = \Z[M_T]\).
We use \(\oplus\) to denote the formal sums in \(K_T\).
By a \emph{rational \(T\)-representation}, we mean a representation of a finite cover of \(T\).
Every rational \(T\)-representation has a class in \(K_T\).
In particular, for \(m, n \in \Q\), we have classes \(\chi(m,n) \in K_T\) of rational characters.
See the discussion before \cite[Proposition~5.3.4]{kre:99} for the need to accommodate rational representations.

\begin{proposition}
  \label{prop:p1uO1}
  Fix \(u \in A  \subset \P U\).
  Suppose \(i\) realises the minimum \(\min_i(\frac{1}{d_i}(r_i^u+b_i))\).
  Then the map
  \begin{equation}\label{eqn:eb}
    \iota^{*} \O(-d_i) \to W_i \otimes \O_{\mathscr M^{\rm res}}
  \end{equation}    
  is non-zero at \(p_1^u\) and its image is spanned by
  \( v_1^{a_i-b_i-r_i}v_2^{r_i} \otimes (v_1 \wedge v_2)^{b_i}\).
\end{proposition}
\begin{proof}
  In the local coordinates introduced in \Cref{sec:fp}, the point \(p_1^u\) lies in the chart
  \[\{(x,y,z,[X:1]) \mid x = yX\} \cong \spec \k[X,y,z]\]
  of the blow-up \(\Bl_{K_u}M\).
  The stacky blow-up \(\mathscr M^{\rm res}\) of this chart is defined by the weighted monomials \((y^{r_i^u+b_i},d_i)\).
  Suppose the minimum in the statement is \(c/d\), where \(\gcd(c,d) = 1\).
  From \eqref{eqn:chart} and \eqref{eqn:blowupmap}, we get the following local chart of \(\mathscr M^{\rm res} \to \Bl_{K_u}M\) at \(p_1^u\):
  \[ [\spec \k[u_1,u_2,u_3] / \mu_d] \to \spec \k[X,y,z],\]
  where the map is defined by
  \[ X \mapsto u_1, \quad y \mapsto u^{d}_2, \quad z \mapsto u_3.\]
  Let \(r = r_i^u\) and \(b = b_i\).
  From the proof of \Cref{cor:resolution}, we know that on \(\spec \k[u_1,u_2,u_3]\), the map \eqref{eqn:eb} is \(y^{-r-b}\) times the original map \(e\) studied in \Cref{prop:zerolocus}.
  From \eqref{eqn:pex} and \eqref{eqn:pexdet} , we see that the map \(e\) is given by
  \[ y^{b}\left((v_1+zv_2)^{a-b-r}y^{r}(Xv_{1}+(1+Xz)v_{2})^{r}   + \cdots \right) \otimes (v_1 \wedge v_2)^b.\]
  Multiplying by \(y^{-r-b}\) and setting \(u_1 = u_2 = u_3 = 0\) yields the result.    
\end{proof}

Let \(\O(-1)_{p_1^u}\) be the class in \(K_T\) of the fiber of \(\iota^{*}\O(-1)\) at \(p_1^u\).
Similarly, let \(N_{p_1^u}\) be the class in \(K_T\) of the normal bundle of \(p_1^u\) in \(\mathscr M^{\rm can}\).
Let \(r^u = \min_i(r_i^u+b_i)/d_i\).
\begin{proposition}\label{prop:p1uall}
  With the notation above, we have
  \begin{align*}
    \O(-1)_{p_1^u} &= \chi(1- r^u, r^u), \text{ and }\\
    N_{p_1^u} &=  \chi(1,-1) \oplus \chi(-1,1) \oplus \chi(-1,1).             
  \end{align*}
\end{proposition}
\begin{proof}
  Suppose \(i\) realises the minimum \(\min_i(r_i^u+b_i)/d_i\).
  \Cref{prop:p1uO1} identifies the fiber of \(\iota^{*}\O(-d_i)\) at \(p_1^u\) with the span of \(v_1^{a_i-r_i}v_2^{r_i} \otimes (v_1\wedge v_2)^{b_i}\), on which \(T\) acts by weights \(a_i-r_i\) and \(b_i+r_i\).
  Dividing through by \(d_i\) yields the first equality.

  The map \(\mathscr M^{\rm can} \to \Bl_AM\) is an isomorphism near \(p_1^u\).
  The normal space at \(p_1^u\) is spanned by \(\frac{\partial}{\partial X}\), \(\frac{\partial}{\partial y}\), and  \(\frac{\partial}{\partial z}\), on which the \(T\) acts by weights \((1,-1)\), \((-1,1)\), and \((-1,1)\), respectively.

\end{proof}

\begin{proposition}\label{prop:p1juO1}
  Fix \(u \in A \subset \P U\).
  Let \(\lambda\) be a vertex of \(\Lambda^u\), say \(\lambda = \frac{1}{d_i}(r_i^u+b_i,b_i)\).
  Then the map
  \begin{equation}\label{eqn:eb2}
    \iota^{*} \O(-d_i) \to W_i \otimes \O_{\mathscr M^{\rm res}}
  \end{equation}
  is non-zero at \(p_{1,j}^u\) and its image is spanned by \(v_1^{a_i-b_i} \otimes (v_1\wedge v_2)^{b_i}\).
\end{proposition}
\begin{proof}
  The proof is parallel to the proof of \Cref{prop:p1uO1}.
  In the local coordinates introduced in \Cref{sec:fp}, consider the chart of \(\Bl_{K_u}M\) given by
  \[\{(x,y,z,[1:Y]) \mid y = xY\} \cong \spec \k[x,Y,z].\]
  Consider the chart of \(\mathscr M^{\rm res} \to \Bl_AM\) given by \eqref{eqn:chart} and \eqref{eqn:blowupmap}:
  \[[\k[u_1,u_2,u_3] / \mu] \to \k[x,Y,z].\]
  Let \(r = r_i^u\) and \(b = b_i\).
  From the proof of \Cref{cor:resolution}, we know that on \(\spec \k[u_1,u_2,u_3]\), the map \eqref{eqn:eb2} is \(x^{-r-b}Y^{-b}\) times the original map \(e\) studied in \Cref{prop:zerolocus}.
  From \eqref{eqn:pex} and \eqref{eqn:pexdet}, we see that the map \(e\) is given by
  \[
   x^bY^b((v_1+zv_2)^{a-b-r}x^r(v_1+(1+Yz)v_2)^{r} + \cdots) \otimes (v_1 \wedge v_2)^b.
  \]
Multiplying by \(x^{-r-b}Y^{-b}\) and setting \(u_1 = u_2 = u_3 = 0\) yields the result.
\end{proof}

% \( v_1^{a_i-b_i} \otimes (v_1 \wedge v_2)^{b_i}\).
Let \(\O(-1)_{p_{\ell,j}^u}\) be the class in \(K_T\) of the fiber of \(\iota^{*}\O(-1)\) at \(p_{1,j}^u\).
Similarly, let \(N_{p_{1,j}^u}\)  be the class in \(K_T\) of the normal bundle of \(p_{1,j}^u\) in \(\mathscr M^{\rm can}\).
Let \(\eta\) and \(\zeta\) be the shortest integral normal vectors to the two rays of \(\Lambda^u\) at the vertex \(\lambda(j)\).
Let \(N = \det(\eta,\zeta)\), so that \(|N|\) is the index of the sub-lattice \(\langle \eta, \zeta \rangle \subset \Z^2\).
\begin{proposition}\label{prop:p1juall}
  With the notation above, we have
  \begin{align*}
    \O(-1)_{p_{1,j}^u} &= \chi(1-\lambda(j)_2,\lambda(j)_2), \text { and }\\
    N_{p_{1,j}^u} &= \chi(\zeta_1/N,-\zeta_1/N) \oplus \chi(-\eta_1/N,\eta_1/N) \oplus \chi(-1,1).
  \end{align*}
\end{proposition}
\begin{proof}
  We use the notation in \Cref{prop:p1juO1}.
  In particular, we let \(i\) be such that \(\lambda(j) = \frac{1}{d_i}(a_i+r_i^u,b_i)\).
  \Cref{prop:p1juO1} shows that \(\O(-1)_{p_{1,j}^u} = \frac{1}{d_i}\chi(a_i,b_i)\), by the same argument as \Cref{prop:p1uall}.
  Since \(d_i = a_i+b_i\), we can re-write this as \(\chi(1-\lambda(j)_2, \lambda(j)_2)\).

  For the second equality, we write \(\mathscr M^{\rm can} \to \Bl_AM\) in charts at \(p_{1,j}^u\) using \eqref{eqn:chart} and \eqref{eqn:blowupmap}:
  \[
   [\spec \k[u_1,u_2,u_3]/\mu] \to \spec \k[x,Y,z],
 \]
 where the map is given by
 \[ x \mapsto u_1^{\eta_1}u_2^{\zeta_1}, \quad Y \mapsto u_1^{\eta_2}u_2^{\zeta_2}, \quad z \mapsto u_3.\]
 The torus \(T\) acts on \(x\), \(Y\), and \(z\) by weights \((0,0)\), \((1,-1)\), and \((1,-1)\), respectively.
 It follows that it must act on \(u_1\), \(u_2\), and \(u_3\) by weights \(\frac{1}{N}(-\zeta_1,\zeta_1)\), \(\frac{1}{N}(\eta_1,-\eta_{1})\), and \((1,-1)\), respectively.
 Since the normal space to \(p_{1,j}^u\) is spanned by \(\frac{\partial}{\partial u_1}\), \(\frac{\partial}{\partial u_2}\), and \(\frac{\partial}{\partial u_3}\), the second equality follows.
\end{proof}
\begin{remark}
    In \Cref{prop:p1juall}, suppose \(\lambda = \lambda(0)\) is the bottom right vertex of \(\Lambda^u\).
    Then one of the two rays incident at \(\lambda\) is \(\lambda + \R_{\geq 0}\), so we may choose \(\eta = (0,1)\).
    In that case, we have a trivial summand \(\chi(0,0)\) in \(N_{p_{1,0}^u}\).
    This summand corresponds to the normal direction in \(\mathcal L_1\), which is fixed by \(T\).
  \end{remark}

  We have now described all the ingredients of the localisation formula for the isolated fixed points.
  We now turn to the \(T\)-fixed lines \(\mathscr L_i\).
  Note that the coarse space of \(\mathscr L_i\) is \(\P^1\).
  Use \(K\) to denote the numerical Grothendieck group (two classes considered equal if they have the same Chern character).
  For a finite cover \(\widetilde T \to T\), a \(\widetilde T\)-equivariant bundle on \(\mathscr L_i\) has a class in \(K_T \otimes K(\mathscr L_1)\).
  For \(a \in \Q\), the notation \(\O(a)\) denotes the class in \(K(\mathscr L_i)\) of a line bundle of degree \(a\).
  
  For \(u \in A\), we denote by \(\lambda^u(j)\) the \(j\)-th vertex of \(\Lambda^u \subset \R^2\) with the convention that the vertices are arranged from the bottom-right to the top-left (in the increasing order by the second coordinate).
  Recall that \(b = \min_i(b_i/d_i)\) and \(r_{\rm gen}^u = \lambda^u(0)_1-b\) and \(r_{\rm gen} = \sum_{u \in A} r_{\rm gen}^u\).
  \begin{proposition}\label{prop:universalsub}
    The class of \(\iota^{*} \O(-1)\) restricted to \(\mathscr L_1\) in \(K_T \otimes K(\mathscr L_1)\) is given by
    \[\iota^{*} \O(-1)|_{\mathscr L_1} = \chi(1-b,b) \otimes \O(-1+2b+r_{\rm gen}).\]
  \end{proposition}
  \begin{proof}
    Recall that the points \(p_{1,0}^u\) lie on \(\mathscr L_1\).
    \Cref{prop:p1juall} applied to \(j = 0\) shows that the fiber of \(\iota^{*}\O(-1)\) at \(p_{1,0}^u\), as a rational \(T\)-representation, is \(\chi(1-b,b)\).
    So, the class of \(\iota^{*}\O(-1)\) restricted to \(\mathscr L_{1}\) is \(\chi(1-b,b) \otimes \O(a)\) for some \(a \in \Q\), which is simply the degree of \(\iota^{*}\O(-1)\) on \(\mathscr L_1\).

    Let \(\pi \colon \mathscr M^{\rm res} \to M\) be the natural map.
    Then
    \[ \iota^{*} \O(-1) = \pi^{*} \O_M(-1) \otimes \O(E),\]
    where \(E \subset \mathscr M^{\rm res}\) is an effective divisor.
    The divisor is characterised by the property that in a neighbourhood of a point \(p \in \mathscr M^{\rm res}\) at which the map \(\iota^{*}\O(-d_i) \to W_i \otimes \O_{\mathscr M^{\rm res}}\) is non-zero, the divisor \(d_i E\) is the vanishing locus of \(e \colon \pi^{*} \O_M(-d_i) \to W_i \otimes \O_{\mathscr M^{\rm res}}\).
    We use this characterisation at \(p_{1,0}^u\) for every \(u \in A \subset \P U\).
    Given \(u \in A\), let \(i\) be such that \(\lambda^u(0) = \frac{1}{d_i}(r_i^u+b_i,b_i)\).
    By \Cref{prop:p1juO1}, the map \(\iota^{*}\O(-d_i) \to W_i \otimes \O_{\mathscr M^{\rm res}}\) is non-zero.
    The proof of \Cref{prop:p1juO1} shows that the vanishing locus of \(e \colon \pi^{*} \O_M(-d_i) \to W_i \otimes \O_{\mathscr M^{\rm res}}\) is cut out in the local coordinates by \(x^{r_i^u+b_i}Y^{b_i} = x^{r_i^u}y^{b_i}\).
    The divisor cut out by \(y\) is the pre-image of the determinant \(\Delta \subset M\).
    The divisor cut out by \(x\) is the pre-image of the exceptional divisor \(E_u \subset \Bl_AM\).
    Therefore, in a neighbourhood of \(p_{1,0}^u\), we have
    \begin{align*}
      E &= \frac{b_i}{d_i} \Delta + \frac{r_i^u}{d_i} E_u \\
        &= b \Delta + r_{\rm gen}^u E_u \
    \end{align*}
    Considering all \(u \in A \subset \P U\), we see that in a neighbourhood of \(\mathscr L_1\), we have
    \[ E = b \Delta + \sum_{u \in A} r_{\rm gen}^u E_u.\]
    On \(\mathscr L_1\), the degree of \(\Delta\) is \(2\) and the degree of \(E_u\) is \(1\).
    The result follows.
  \end{proof}

Let \(N_1 \in K_T \otimes K(\mathscr L_1)\) be the class of the normal bundle of \(\mathscr L_1^{\rm can} \subset \mathscr M^{\rm can}\).
  For \(u \in A \subset \P U\), if \(\Lambda^u\) has at least two vertices, set
  \[ s^u = 1-\frac{\lambda^u(0)_1 - \lambda^u(1)_1}{\lambda^u(0)_2 - \lambda^u(1)_2}.\]
  Otherwise, set \(s^u = 1\).
  Let \(s = \sum_{u \in A} s^u\).
  \begin{proposition}\label{prop:normal}
    With the notation above, the class of the normal bundle \(N_1\) is equal to
    \[ \chi(-1,1) \otimes \left( \O \oplus \O(2-s))\right).\]
  \end{proposition}
  \begin{proof}
    For simplicity, we drop the superscript ``\(\rm can\)'', but alert the reader that it is important that we are working with \(\mathcal M^{\rm can}\) and not \(\mathcal M^{\rm res}\).
    We use the local coordinates introduced in \Cref{sec:fp}.
    At a generic point \((x,0,0) \in \mathscr L_1\), the normal bundle is spanned by \(\frac{\partial}{\partial y}\) and \(\frac{\partial}{\partial z}\) on which \(T\) acts by weights \((-1,1)\).
    Therefore, the class of \(N_1\) is \(\chi(-1,1)\) times the class in \(K(\mathscr L_1)\) of the normal bundle.
    We simply need to find the degree of \(N_1\).
    
    Let \(\pi \colon \mathscr M \to \Bl_AM\) be the natural map.
    Let \(\widetilde L_1 \subset \Bl_AM\) be the proper transform of \(L_1\).
    Consider the sequence
    \[ 0 \to N_{1} \xrightarrow{d \pi} \pi^{*}N_{\widetilde L_1/\Bl_A M} \to Q \to 0,\]
    so that \(Q\) is supported on \(\{p_{1,0}^u \mid u \in A\}\).
    Let \(\eta = (0,1)\) and \(\zeta\) be the shortest integer normal vectors to the two rays of \(\Lambda^u\) at \(\lambda(0)\).
    Then, in a neighbourhood of \(p_{1,0}^u\), the map \(\pi\) is
    \[
   [\spec \k[u_1,u_2,u_3] / \mu] \to \spec \k[x,Y,z],
 \]
 where
 \[ x \mapsto u_2^{\zeta_1}, \quad Y \mapsto u_1u_2^{\zeta_2}, \quad z \mapsto u_3, \]
 and \(\mu\) is a cyclic group of order \(\zeta_1\).
 At \(p_{1,0}^u\), the dual of \(N_{\widetilde L_1/\Bl_A M}\) is spanned by \(dY\) and \(dz\), whereas the dual of \(N_{1}\) is spanned by \(du_1\) and \(du_3\).
 On \(\mathscr L_1\), which is cut out by \(u_1 = u_3 = 0\), we have
  \[ dY = u_2^{\zeta_2} du_1 \text{ and } dz = du_3.\]
  Therefore, the pull-back of \(Q\) to \(\k[u_1,u_2,u_3]\) has length \(\zeta_2\).
  But since the order of \(\mu\) is \(\zeta_1\), the degree of \(Q\) at \(p_{1,0}^u\) is \(\zeta_2/\zeta_1\).
  If \(\Lambda^u\) has only one vertex, then \(\zeta = (1,0)\), so \(\zeta_2 = 0\).
  Otherwise, \(\zeta_2/\zeta_1\) is the negative of the reciprocal of the slope of the line joining \(\lambda(0)\) and \(\lambda(1)\), which is precisely \(s^u-1\).
  Therefore, we conclude that
  \[ \deg N_{1} = \deg N_{\widetilde L_1 / \Bl_A M} - \sum_{u \in A} (s^u-1).\]
  But we also know that
  \[ \deg N_{\widetilde L_1 / \Bl_A M} = \deg N_{L_1 / M} - \sum_{u \in A} 1 = 2 - \sum_{u \in A} 1.\]
  Combining the two yields the proposition.
\end{proof}
\subsection{Proof of the main theorem}
We now have the tools to prove \Cref{thm:main}.
Let \(T \subset \GL(V)\) be a maximal torus.
Let \(N = \dim W\).
By \Cref{prop:weighted-complete-param}, we have
\[
  |\Gamma| \cdot [\Orb(w)] = \int_{\mathscr M^{\rm res}} \frac{c_{N}(W)}{\iota^{*}c_1 \O_{\sP W}(-1)}.
\]
The pull-back along \(\mathscr M^{\rm res} \to \mathscr M^{\rm can}\) identifies the rational Chow groups, and the push-forward of the fundamental class of \(\mathscr M^{\rm res}\) is the fundamental class of \(\mathscr M^{\rm can}\).
Therefore, we may replace \(\mathscr M^{\rm res}\) in the integral by \(\mathscr M^{\rm can}\).

From \Cref{sec:fp}, recall that the \(T\)-fixed points of \(\mathcal M^{\rm can}\) consist of the line \(\mathscr L_1^{\rm can}\), the points \(p_1^u\) for \(u \in A\), the points \(p_{1,j}^u\) for \(u \in A\) and \(j = 1, \dots, k^u\), where \(k^u\) is the number of vertices of the Newton polygon \(\Lambda^u\), and their analogues where the subscript \(1\) is replaced by \(2\).
For \(\ell = 1,2\), we denote by \(N_{\ell}\) the normal bundle of \(\mathscr L_{\ell}^{\rm can}\), by \(N_{p_{\ell}^u}\) the normal space of \(p_{\ell}^u\), and by \(N_{p_{\ell,j}^u}\) the normal space of \(p_{\ell,j}^u\).
Let us write \( \xi = c_{N}(W) / \iota^{*}c_1 \O_{\sP W}(-1)\).
By the localisation formula, we have the equality of \(T\)-equivariant classes
\begin{equation}\label{eqn:bigsum}
  \begin{split}
    \int_{\mathscr M^{\rm can}} \xi = \int_{\mathscr L^{\rm can}_1} \frac{\xi}{c_{2}(N_1)} + \sum_{u \in A}\int_{p_1^u} \frac{\xi}{c_3(N_{p_1^u})} + \sum_{u \in A} \sum_{j = 1}^{k^u} \int_{p_{1,j}^u} \frac{\xi}{c_3(N_{p_{1,j}^u})} + \cdots
  \end{split}
\end{equation}
where \(\cdots\) denotes the sum of analogous integrals over \(\mathscr L_2^{\rm can}\) and \(p_2^u\) and \(p_{2,j}^u\).

Let us now evaluate each term in \eqref{eqn:bigsum}, starting with the integral over \(\mathscr L_1^{\rm can}\).
Denote by \(h \in A^1(\P^1)\) the class of a point.
Using \Cref{prop:universalsub} and \Cref{prop:normal} (and the notation there), we have
\begin{align*}
  \frac{1}{c_N(W)}&\int_{\mathscr L_1^{\rm can}} \frac{\xi}{c_2(N_1)} = \int \frac{1}{c_1(\O(-1)) \cdot c_2(N_1)} \\
  &= \int ((1-b)v_1+bv_2+ (2b+r_{\rm gen}-1)h)^{-1}(v_2-v_1)^{-1}(v_2-v_{1}+(2-s)h)^{-1}.
\end{align*}
The integral is the coefficient of \(h\) in the expansion of the integrand as a power series in \(h\).
To find it, we formally differentiate with respect to \(h\) and set \(h = 0\) to obtain
\begin{equation}\label{eqn:Lcontrib}
  \begin{split}
    \frac{1}{c_N(W)}\int_{\mathscr L_1^{\rm can}} \frac{\xi}{c_2(N_1)}&=
                                                                                                      (2-s)((1-b)v_1+bv_2)^{-1}(v_1-v_2)^{-3} \\ & \qquad - (2b+r_{\rm gen}-1)((1-b)v_1+bv_2)^{-2}(v_1-v_2)^{-2}.
  \end{split}
\end{equation}
The analogous integral over \(\mathscr L_2^{\rm can}\) is obtained by switching \(v_1\) and \(v_2\).

Let us turn to the integral over \(p_1^u\).
By \Cref{prop:p1uall} (and the notation there), we have
\begin{equation}\label{eqn:p1ucontrib}
  \begin{split}
    \frac{1}{c_N(W)}\int_{p_1^u} \frac{\xi}{c_3(N_1)} &=
                                                                                        \int {((1-r^u) v_1 + r^u v_2)^{-1}(v_1-v_2)^{-3}} \\
                                                                                      &= {((1-r^u) v_1 + r^u v_2)^{-1}(v_1-v_2)^{-3}}
  \end{split}
\end{equation}
The analogous integral over \(p_2^u\) is obtained by switching \(v_1\) and \(v_2\).

Finally, let us compute the integral over \(p_{1,j}^u\).
By \Cref{prop:p1juall} (and the notation there), we have
\begin{equation}\label{eqn:p1jucontrib}
  \begin{split}
    \frac{1}{c_N(W)}\int_{p_{1,j}^u} \frac{\xi}{c_3(N_1)} &=
                                                                                        \int {N^2}{((1-\lambda(j)_2) v_1 + \lambda(j)_2 v_2)^{-1}(v_1-v_2)^{-3}\zeta^{-1}_1\eta^{-1}_1} \\
                                                                                      &= {|N|}{\zeta^{-1}_1\eta^{-1}_1((1-\lambda(j)_2) v_1 + \lambda(j)_2 v_2)^{-1}(v_1-v_2)^{-3}}.
  \end{split}
\end{equation}
In the last equality, we have used that \(p_{1,j}^u \in \mathscr M^{\rm can}\) has a stabiliser of order \(|N|\), and hence the integral divides the integrand by \(|N|\).
The analogous integral over \(p_{2,j}^u\) is obtained by switching \(v_1\) and \(v_2\).

The expression in \Cref{thm:main} is the sum of the contributions from \eqref{eqn:Lcontrib}, \eqref{eqn:p1ucontrib}, \eqref{eqn:p1jucontrib}, and their analogues with \(v_1\) and \(v_2\) switched.

\section{Applications}\label{sec:applications}
\subsection{Orbits of elliptic fibrations}\label{sec:ellipticfibrations}
Recall that an element \((A,B) \in \Sym^{4n}(V) \oplus \Sym^{6n}(V)\) gives rise to an elliptic fibration
\[ \pi \colon E \to \P^1\]
defined locally by the Weierstrass equation
\[ y^2 = x^3 + Ax + B.\]
Given \(u \in \P^{1}\), recall that \(r_1^u\) is the order of vanishing of \(A\) at \(u\) and \(r_2^u\) is the order of vanishing of \(B\) at \(u\).
We are now ready to prove \Cref{thm:ellipticfibrations}, which computes the degree of the orbit closure of \((A,B)\).
\begin{proof}[Proof of \Cref{thm:ellipticfibrations}]
  Let \(w = (A,B) \in W = \Sym^{4n}(V) \oplus \Sym^{6n}(V)\) be non-zero.
  In the notation of \Cref{thm:main}, we have \(b = 0\).
  For every \(u \in \P^1\), the Newton polygon \(\Lambda^u\) has only one possible shape.
  It is a translated quadrant \(\lambda + \mathbf{R}_{\geq 0}\) whose vertex \(\lambda\) is
  \[ \lambda = \left(\min\left(\frac{1}{4n}r_1^u, \frac{1}{6n}r_2^u\right), 0\right) = \left(\frac{c(u)}{2n},0\right).\]
  In the notation of \Cref{thm:main}, we have \(r^u = r_{\rm gen}^u = c(u)/2n\) and \(s^u = 1\).

  Note that \(\P W\) is the quotient of  \(W - 0\) by the \(\Gm\) acting by weights \(2\) and \(3\).
The central \(\Gm \subset \GL V\) acts by weights \(4n\) and \(6n\).
Therefore, the equivariant class for the first \(\Gm\) is obtained from the \(\GL V\)-equivariant class by the specialisation \(v_1 = v_2 = \frac{h}{2n}\).
With these substitutions, \Cref{thm:ellipticfibrations} follows from \Cref{thm:main}.
\end{proof}

\subsection{Orbits of rational self maps}\label{sec:ratmaps}
Recall that elements in a Zariski open subset of \(\Hom(V, \Sym^nV)\) give rise to maps \(f \colon \P V \to \P V\) of degree \(n\).
We have an isomorphism of \(\GL V\)-representations
\begin{equation}\label{eqn:splitting}
  \Hom(V, \Sym^n V) = \Sym^{n-1} V \oplus \Sym^{n+1} V \otimes \det V^{-1}.
\end{equation}
The first projection \(\Hom(V, \Sym^nV) \to \Sym^{n-1} V\) is the contraction.
The second projection arises as the composite
\[\Hom(V, \Sym^nV) \otimes \det V = V \otimes \Sym^nV \to \Sym^{n+1}V\]
where the first map arises from the isomorphism \(V^{*} \otimes \det V = V\) and the second map is the multiplication.
The element in \(\Sym^{n+1}V\) in the second projection defines the scheme theoretic fixed locus of \(f\).
I do not know a similar geometric interpretation of the first projection.

Fix a basis \(x,y\) of \(V\) with the dual basis \(x^{*}, y^{*}\) of \(V^{*}\).
\begin{proposition}\label{prop:fgab}
  Let \(f = x^{*} \otimes F(x,y) + y^{*} \otimes G(x,y)\), where \(F, G \in \Sym^nV\) are polynomials of degree \(n\).
  Let \(f\) correspond to \((I,J \otimes x^{*} \wedge y^{*}) \in \Sym^{n-1}V \otimes \Sym^nV \otimes \det V^{-1}\).
  Then, up to non-zero scalar multiples, we have
  \[I = \frac{\partial}{\partial x}F(x,y) + \frac{\partial}{\partial y}G(x,y) \text{ and }  J = yF(x,y) - xG(x,y).\]
  In the other direction, we have
  \[ \frac{1}{n+1} F = \frac{\partial J}{\partial y} + xI \text{ and } \frac{1}{n+1} G = yI - \frac{\partial J}{\partial x}.\]
\end{proposition}
\begin{proof}
  It is enough to check that the construction of \(I\) and \(J\) is \(\GL(2)\)-equivariant.
  We leave this to the reader.
  The other direction follows from the first using Euler's formula
  \[ x\frac{\partial *}{\partial x} + y\frac{\partial *}{\partial y}  = \deg(*) \cdot *.\]
\end{proof}

Fix a non-zero \(\delta \in \det V\).
\begin{proposition}\label{prop:ab}
  Suppose \(f \in \Hom(V, \Sym^nV)\) defines a rational map \(\P V \to \P V\) of degree \(n\) and corresponds to \((I, J \otimes \delta^{-1})\) under an isomorphism \eqref{eqn:splitting}.
  If \(J\) vanishes to order at least \(2\) at \(u \in \P V\), then \(I\) does not vanish at \(u\).
\end{proposition}
\begin{proof}
  Write \(f = x^{*} \otimes F(x,y) + y^{*} \otimes G(x,y)\) in coordinates.
  Since \(f\) defines a map of degree \(n\), the polynomials \(F\) and \(G\) have no common factor.
  If \(J\) vanishes to order at least 2 at \(u\), then both partials of \(J\) vanish to order at least 1 at \(u\).
  Since at least one of \(F\) or \(G\) does not vanish at \(u\), we see from \Cref{prop:fgab} that \(I\) cannot vanish at \(u\).
\end{proof}

We are now ready to prove \Cref{thm:ratmaps}, which computes the equivariant orbit class of a rational map.
\begin{proof}[Proof of \Cref{thm:ratmaps}]
  Let \(f \in \Hom(V, \Sym^nV)\) define a rational map \(\P V \to \P V\) of degree \(n\).
  Let \(f\) correspond to \((I, J \otimes \delta^{-1})\) under an isomorphism \eqref{eqn:splitting}.
  We apply \Cref{thm:main}, taking \(A = V(J)\) to be the set of fixed points of \(\P V \to \P V\).
  We have \(b = -1/(n-1)\).
  By \Cref{prop:ab}, for \(u \in A\), the Newton polygon \(\Lambda^u\) has two possible shapes (see \Cref{fig:twoshapes}).
  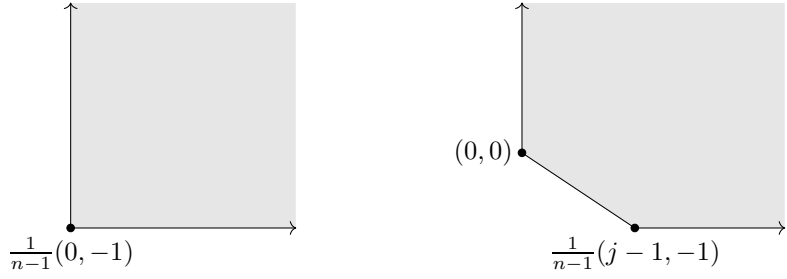
\begin{figure}[ht]
    \centering
    \begin{tikzpicture}
      \begin{scope}
        \draw[draw=white, fill=black!10!white] (0,-1) -- (0,2) -- (3,2) -- (3,-1) -- (0,-1);
        \draw[fill] (0,-1) circle (0.05) node[below] {\(\frac{1}{n-1}(0,-1)\)};
        \draw[->] (0,-1) -- (0,2);
        \draw[->] (0,-1) -- (3,-1);
      \end{scope}
      \begin{scope}[xshift=6cm]
        \draw[draw=white, fill=black!10!white] (0,0) -- (1.5,-1) -- (3.5,-1) -- (3.5,2) -- (0,2) -- (0,0);
        \draw[fill] (0,0) circle (0.05) node[left] {\((0,0)\)};
        \draw[fill] (1.5,-1) circle (0.05) node[below] {\(\frac{1}{n-1}(j-1,-1)\)};
        \draw (0,0) -- (1.5,-1);
        \draw[->] (0,0) -- (0,2);
        \draw[->] (1.5,-1) -- (3.5,-1);
      \end{scope}
    \end{tikzpicture}
    \caption{\(\Lambda^u\) for a simple fixed point \(u\) (left) and a fixed point of order \(j \geq 2\) (right)}
    \label{fig:twoshapes}
  \end{figure}
  If \(u\) is a simple fixed point, then its only vertex is
  \[\lambda^u(0) = \frac{1}{d_2} (r_2^u+b_2,b_2) = \frac{1}{n-1}(0,-1).\]
  In this case, \(r_{\rm gen}^u = 1/(n-1)\) and \(r^u = 0\) and \(s^u = 1\).
  If \(u\) is a fixed point of order \(j^u \geq 2\), then the vertices of \(\Lambda^u\) are
  \[
    \begin{split}
      \lambda^u(0) &= \frac{1}{d_2}(r_2^u+b_2,b_2) = \frac{1}{n-1}{(j^u-1,-1)} \text{ and }\\
      \lambda^u(1) &= \frac{1}{d_1}(b_1,b_1) = (0,0).
    \end{split}
  \]
  In this case, \(r_{\rm gen}^u = j^u/(n-1)\) and \(r^u = 0\) and \(s^u = j^u\).
  In the notation of \Cref{thm:main}, we have
  \begin{align*}
    F &= 2(n-1)(n v_1 - v_2)^{-1}(v_1-v_2)^{-3} + (n+1)(n v_1 - v_2)^{-2}(v_1-v_2)^{-2}, \text{ and } \\
    G^u &= v_1^{-1}(v_1-v_2)^{-3}-j^{u}(n-1)(n v_1 - v_2)^{-1}(v_1-v_2)^{-3}- j^{u}(n v_1 - v_2)^{-2}(v_1-v_2)^{-2},
  \end{align*}
  and for a higher order fixed point \(u\),
  \[
   H^u(1) = (j^u-1) v_1^{-1}(v_1-v_2)^{-3}.
 \]
 Summing up and multiplying by \(c_N(W)\) yields the class
 \[ n(n+1)(n-1)^2\prod_{j = 1}^{n-2}(jv_1+(n-1-j)v_2) \prod_{j = 1}^n((j-1) v_1 + (n-j)v_2). \]
 Note that \(\P \Hom(V, \Sym^nV)\) is the quotient of \(\Hom(V, \Sym^nV) - 0\) by \(\Gm\) acting by weight one.
 The central \(\Gm \subset \GL V\) acts by weight \(n-1\).
 Therefore, the weight one \(\Gm\) equivariant class is obtained from the \(\GL V\)-equivariant class by specialising to \(v_1 = v_2 = 1/(n-1)\).
 The stabiliser group \(\Gamma \subset \GL V\) in \Cref{thm:main} and the stabiliser group \(\overline \Gamma \subset \PGL V\) in \Cref{thm:ratmaps} are related by the exact sequence
 \[ 1 \to \mu_{n-1} \to \Gamma \to \overline \Gamma \to 1.\]
 So we must divide the class given by \Cref{thm:main} by \((n-1)\).
 Specialising to \(v_1 = v_2 = 1/(n-1)\) and dividing by \((n-1)\) gives \(n(n+1)(n-1)\).
 \end{proof}

 \appendix
 \section{Equivariant classes of torus orbits}\label{sec:torus}
 Fix an algebraic torus \(T = \Gm^d\).
 We compute \(T\)-equivariant fundamental classes of orbits in \(T\)-representations.
 Let \(M = \Hom(T, \Gm)\) be the character group of \(M\) and set \(N = \Hom(M, \Z)\).
 Identify \(A_T = \Sym(M_{\Q})\).
 
 Fix a \(T\)-representation \(W\) and a \(w \in W\).
 Write \(W = \bigoplus_{i = 1}^n W_i\), where \(W_i\) is one-dimensional on which \(T\) acts by the character \(\chi_i \in M\).
Fix a non-zero \(w = (w_1, \dots, w_n) \in W\).

 We recall the notion of equivariant multiplicity from \cite{bri:97}.
 Given a polyhedral rational pointed cone \(\sigma \subset M_{\R}\), denote by \(\sigma^\vee \subset N_{\R}\) the dual cone.
 Since \(\sigma\) is pointed, \(\sigma^{\vee}\) has non-empty interior.
 Given \(\lambda\) in the interior of \(\sigma^{\vee}\), let \(P_{\sigma}(\lambda)\) be the convex polytope
 \[ P_{\sigma}(\lambda) = \{x \in \sigma \mid \langle  x, \lambda \rangle \leq 1\}.\]
 There exists a unique rational function \(e_{\sigma} \in \operatorname{frac}{\Sym(M}_{\Q})\) such that for every \(\lambda\) in the interior of \(\sigma^{\vee}\), we have
 \[ e_{\sigma}(\lambda) = d! \cdot \operatorname{Vol} P_{\sigma}(\lambda).\]
The function \(e_{\sigma}\) is called the \emph{equivariant multiplicity} associated to \(\sigma\) (see \cite[\S~5.2]{bri:97}).

 Let \(\sigma \subset M_{\R}\) be the closed convex cone spanned by \(\{\chi_i \mid w_i \neq 0\}\).
 Let \(\Orb(w) \subset W\) be the closure of the \(T\)-orbit of \(w\) and \([\Orb(w)]\) its fundamental class in \(A_{T}(W) = A_{T}\).
 \begin{theorem}\label{thm:torus}
   In the setup above, if \(\sigma\) contains a line, then \([\Orb(w)] = 0\).
   Otherwise,
   \[ [\Orb(w)] = e_{\sigma} \cdot c_n(W).\]
 \end{theorem}
 \begin{proof}
   If \(\sigma\) contains a line, then \(0 \in W\) is not in \(\Orb(w)\).
   As a result, the pull-back of \([\Orb(w)]\) to \(A_{T}(0)\) vanishes.
   But the pull-back \(A_T(W) \to A_T(0)\) is an isomorphism, so \([\Orb(w)]\) vanishes.

   Assume that \(\sigma\) contains no line, that is, it is pointed.
   Let \(X\) be the affine toric variety
   \[ X = \spec \k[M \cap \sigma].\]
   It is easy to check that the map \(T \to W\) that sends \(t \in T\) to \(t \cdot w\) extends to a proper morphism \(i \colon X \to W\).
   Then \([\Orb(w)] = i_{*}[X]\).

   To compute the push-forward, we use localisation \cite[\S~4.2 Corollary]{bri:97}.
   Let \(0_X \in X\) and \(0_W \in W\) be the origins.
   Then we have
   \[ [X] = e_{\sigma} \cdot [0_X],\]
   and hence
   \[ i_{*}[X] = e_{\sigma} \cdot [0_W].\]
   Since \([0_{W}] = c_{n}(W)[W]\in A_{T}(W)\), the theorem follows.
 \end{proof}

 \begin{example}
   Let \(T = \G_m^3\) act on \(W = \C^{4}\) by the characters \((0,0,1), (0,1,1), (1,0,1),\) and \((1,1,1)\).
   Take \(w = (1,1,1,1)\).
   Let \(x,y,z\) be the standard basis vectors of \(M=\Hom(T,\G_m)\).
   Let \(\sigma \subset M_{\R}\) be the cone spanned by the four characters.
   Given \(\lambda = (a,b,c) \in \sigma^{\vee} \subset N_{\R}\), we compute
   \[ 3! \cdot \operatorname{vol}(P_{\sigma}(a,b,c)) = \frac{1}{c(b+c)(a+b+c)} + \frac{1}{c(a+c)(a+b+c)}.\]
   Since \(a = \langle x, \lambda \rangle\) and \(b = \langle  y, \lambda \rangle\) and \(c = \langle z, \lambda \rangle\), the equivariant multiplicity function is
   \begin{align*}
     e_{\sigma} &= \frac{1}{z(y+z)(x+y+z)} + \frac{1}{z(x+z)(x+y+z)}\\
&= \frac{x+y+2z}{z(x+z)(y+z)(x+y+z)}.                  
   \end{align*}
   By \Cref{thm:torus}, the equivariant fundamental class of \(\Orb(w)\) is \(x+y+2z\).
   Indeed, in this case, \(\Orb(w) \subset W\) is the quadric hypersurface cut out by \(w_1w_4-w_2w_3\), a polynomial with character \((1,1,2)\).
 \end{example}
 
\bibliography{references}
\bibliographystyle{abbrv}

% For references about

% stacky blow-ups: [[id:aafd57a8-cfc3-4da4-8e01-f752f3973695][Smooth toric Deligne-Mumford stacks]], [[id:eb45da26-691e-4d83-b69f-132b7eda67e0][Weighted blow-ups]]

% stacky Atiyah-Bott: [[id:a0b28674-7376-4e86-93a9-3e2b3b96f999][Localization theorems for algebraic stacks]], [[id:98fa87b3-164f-4c27-88da-8fe059c62cb4][Cycle groups for Artin stacks]]
% 

% You do not want a weighted blow up defined by a Rees algebra.
% You want an iterated weighted blow up.
% The iterated one corresponds to the stacky fan.

% See [[id:209dce03-7a8b-4ea6-9668-b1821a5ca3a6][Logarithmic resolution via multi-weighted blow-ups]]
\end{document}